\documentclass[smallextended]{svjour3} 
\smartqed 
\usepackage{graphicx}
\usepackage{url}
\journalname{COAP}
\usepackage{graphicx,amssymb}
\usepackage{mathtools}
\usepackage[ruled,lined]{algorithm2e}
\usepackage{amsmath}
\usepackage{enumerate}
\usepackage{epstopdf}

\usepackage{marvosym}
\newcommand{\envelope}{(\raisebox{-.5pt}{\scalebox{1.45}{\Letter}}\kern-1.7pt \hspace{0.7mm})}

\newcommand{\prox}{{\text{prox}}}
\newcommand{\sgn}{\text{sgn}}
\newcommand{\supp}{\text{supp}}
\newcommand{\lbr}{\left(}
\newcommand{\rbr}{\right)}
\newcommand{\calP}{\mathcal{P}}
\newcommand{\calQ}{\mathcal{Q}}

\DeclarePairedDelimiter{\ceil}{\lceil}{\rceil}

\newtheorem{assump}{Assumption}

\begin{document}


\title{Local and Global Convergence of a General Inertial Proximal Splitting Scheme\footnote{This manuscript is related to the preprint arxiv:1502.02281, entitled: ``A Lyapunov Analysis of FISTA with Local Linear Convergence for Sparse Optimization". This manuscript is a thoroughly revised and rewritten version which includes several new results.}}


\author{Patrick R. Johnstone \and  Pierre Moulin}

\institute{Patrick R. Johnstone \envelope \and
	Pierre Moulin  \at
             Coordinated Science Laboratory, University of Illinois,
              Urbana, IL 61801, USA\\
              e-mail: prjohns2@illinois.edu, moulin@ifp.uiuc.edu
}

\date{Received: date / Accepted: date}

\maketitle
\begin{abstract}
This paper is concerned with convex composite minimization problems in a Hilbert space. In these problems, the objective is the sum of two closed, proper, and convex functions where one is smooth and the other admits a computationally inexpensive proximal operator. We analyze a general family of inertial proximal splitting algorithms (GIPSA) for solving such problems. We establish finiteness of the sum of squared increments of the iterates and optimality of the accumulation points. Weak convergence of the entire sequence then follows if the minimum is attained.  
Our analysis unifies and extends several previous results.

We then focus on $\ell_1$-regularized optimization, which is the ubiquitous special case where the nonsmooth term is the $\ell_1$-norm. For certain parameter choices, GIPSA is amenable to a local analysis for this problem. For these choices we show that GIPSA achieves finite ``active manifold identification", i.e. convergence in a finite number of iterations to the optimal support and sign, after which GIPSA reduces to minimizing a local smooth  function. Local linear convergence then holds under certain conditions. We determine the rate in terms of the inertia, stepsize, and local curvature. Our local analysis is applicable to certain recent variants of the Fast Iterative Shrinkage-Thresholding Algorithm (FISTA), for which we establish active manifold identification and local linear convergence. Our analysis motivates the use of a momentum restart scheme in these FISTA variants to obtain the optimal local linear convergence rate.

\end{abstract}

\section{Introduction}
The primary problem considered in this paper is to
\begin{eqnarray}
\underset{x\in\mathcal{H}}{\hbox{minimize}}\
F(x)= f(x)+g(x)
\label{prob:1}
\end{eqnarray}
where $\mathcal{H}$ is a Hilbert space over the real numbers, the functions $f,g:\mathcal{H}\to(-\infty,+\infty]$ are proper, convex and closed, and in addition $f$ is 
differentiable everywhere and has a Lipschitz continuous gradient. This problem has come under considerable attention in recent years due to its many applications in areas such as machine learning, compressed sensing and image processing  \cite{bach2011convex,Hale:2008,Beck:2009,prox_signalProcessing,cevher2014convex,tibshirani2013lasso,afonso2010fast,choi2010compressed,lustig2007sparse}. Of particular interest in this paper will be the special case where the nonsmooth term is the $\ell_1$-norm, i.e.
\begin{eqnarray}
\underset{x\in\mathbb{R}^n}{\hbox{minimize}}\
\left\{
f(x)+\rho\|x\|_1
\right\}
\label{prob:sparse}
\end{eqnarray} 
where $\rho> 0$, and $\|x\|_1=\sum_{i=1}^n |x_i|$. As has been widely recognized the $\ell_1$-norm encourages ``sparse" solutions, i.e. solutions with few nonzero elements, which is its primary attraction \cite{bach2011convex,tibshirani2013lasso}. A special case of Prob.~(\ref{prob:sparse}) is
\begin{eqnarray}
\underset{x\in\mathbb{R}^n}{\hbox{minimize}}\
\left\{
\frac{1}{2}\|b-Ax\|_2^2+\rho\|x\|_1
\right\}
\label{prob:l1LS}
\end{eqnarray}
with $A\in\mathbb{R}^{m\times n}$ and $b\in\mathbb{R}^m$, which is often referred to as sparse least-squares, sparse regression, basis pursuit, or lasso and is of vital importance in many areas \cite{kim2007interior,Hale:2008,tibshirani2013lasso,Beck:2009}. Other important instances of Prob.~(\ref{prob:1}) include least-squares with a total-variation \cite{chambolle2015convergence}, nuclear norm \cite{liang2014local}, or group-sparse \cite{bach2011convex} regularizer, and minimization of a convex function constrained to a closed and convex set.

\subsection{Background}
The increasing size of Problems (\ref{prob:1})--(\ref{prob:l1LS}) in modern applications is driving the need for computationally inexpensive and scalable algorithms to find their solutions. In modern applications the number of variables and the number of data can be in the millions \cite{kim2007interior,bach2011convex}. 
\emph{First-order splitting methods} for solving optimization problems including (\ref{prob:1}) are simple and computationally inexpensive, and address the problem by splitting it into simpler subproblems.  
While the overall objective $F$ in Prob.~(\ref{prob:1}) may not have desirable properties, each component of the sum can be handled. The function $f$ is smooth which means it can be processed via its gradient, and many popular nonsmooth regularizers can be processed via a computationally tractable proximal operator \cite{prox_signalProcessing}. Importantly, first-order methods do not rely on or approximate second order information, which may be prohibitively expensive in high dimensions. 
The concept of splitting has also been applied to more complicated objectives  \cite{raguet2013generalized,condat2013primal,eckstein2012augmented}. 
These techniques can also be viewed in the broader context of montone inclusion problems and variational inequalities which includes convex optimization as a special case \cite{raguet2013generalized,condat2013primal,passty1979ergodic,lorenz2014inertial,bauschke2011convex,davis2014convergence}.

The celebrated first-order splitting method for Prob.~(\ref{prob:1}) is the \emph{proximal forward-backward splitting algorithm} (FBS) \cite{passty1979ergodic,combettes2005signal}. 
For this method the convergence rate of the objective function to the optimal value is as good as if the nonsmooth component were not present. Weak convergence of the iterates is also guaranteed and linear convergence occurs on strongly convex problems \cite[Cor.~27.9, Ex.~27.12]{bauschke2011convex}.  
FBS is also commonly referred to as the proximal gradient method. For the special case of Problems (\ref{prob:sparse})--(\ref{prob:l1LS}) it is often referred to as the iterative shrinkage and soft-thresholding algorithm (ISTA) due to the form of the proximal operator with respect to the $\ell_1$-norm.
Other state-of-the-art approaches to solving Prob.~(\ref{prob:1}) and Problems (\ref{prob:sparse})--(\ref{prob:l1LS}) in particular include coordinate descent \cite{friedman2010regularization}, ADMM \cite{afonso2010fast}, and stochastic methods \cite{shalev2011stochastic}. 

Another class of methods of particular interest in this paper are \emph{inertial methods} (a.k.a. momentum methods). These are iterative schemes for solving monotone inclusion and optimization problems, as well as computing fixed points, which often have connections to systems of differential equations (e.g. \cite{alvarez2000minimizing,attouch2014dynamical,su2014differential,PolyakIntro,polyak1964some}). Their defining property is that the next iterate depends on more than one previous iterate (i.e. they are multistep).  
A very early example is due to Polyak \cite{polyak1964some}, who introduced the heavy ball with friction method for minimizing a strongly convex quadratic function which can greatly improve upon the convergence speed of the simple gradient method  (see also \cite[p. 65]{PolyakIntro}). The conjugate gradient method is inertial, as are Nesterov's celebrated accelerated methods, and their variants and extensions \cite{nesterov1983method,Beck:2009,guler1992new,tseng2008accelerated}. Inertial methods typically have the same per-iteration complexity as their noninertial counterparts. However in certain contexts they can be significantly faster   \cite{PolyakIntro,polyak1964some,nesterov2004introductory,Beck:2009}.

\subsection{Contributions}
In this paper we consider the following suite of first-order splitting schemes for solving Prob.~(\ref{prob:1}) which we will refer to as the \emph{General Inertial Proximal Splitting Algorithm} (GIPSA). For all $k\in\mathbb{N}$ compute:
\begin{eqnarray}
\label{IFBS_yupdate}
y^{k+1} &=& x^k + \beta_k(x^k - x^{k-1}),\\
\label{IFBS_zupdate}
z^{k+1} &=& x^k + \alpha_k(x^k - x^{k-1}),\\
\label{IFBS_xupdate}
x^{k+1} &=&\prox_{\lambda_k g}\left(y^{k+1}-\lambda_k\nabla f(z^{k+1})\right) 
\end{eqnarray} 
with $x^{0},x^{1}\in\mathcal{H}$ (typically $x^0=x^1$). The sequences $\{\alpha_k,\beta_k,\lambda_k\}_{k\in\mathbb{N}}$ are in $\mathbb{R}_{+}$ and the proximal operator $\prox_{g}:\mathcal{H}\to\mathcal{H}$ will be properly defined in Sec.~\ref{sec:prox}. Throughout the paper we will refer to $\{\alpha_k,\beta_k\}_{k\in\mathbb{N}}$ as the ``inertia parameters" and $\{\lambda_k\}_{k\in\mathbb{N}}$ as the ``stepsize". FBS (excluding noninertial relaxation parameters) corresponds to GIPSA with $\alpha_k$ and $\beta_k$ set to $0$ for all $k$ (i.e. no inertia). The subclass of algorithms with $\alpha_k = \beta_k$ for all $k$ (i.e. $z^k=y^k$) will be called \emph{Inertial Forward-Backward Splitting} (I-FBS) in this paper, since it corresponds to an inertial version of the Krasnosel'ski\u{i}-Mann (KM) iterations applied to the \emph{Forward-Backward operator} (see Sec.~\ref{sec_assump_optimal} for definition). 

The main motivation for studying (\ref{IFBS_yupdate})--(\ref{IFBS_xupdate}) is that it unifies several existing schemes which correspond to particular parameter choices \cite{alvarez2000minimizing,Hale:2008,lorenz2014inertial,MoudafiOliny,mainge2008convergence,chambolle2015convergence,chambolle2014ergodic}. Thus our global convergence analysis of GIPSA unifies and extends the prior art. Certain special cases of GIPSA (e.g. \cite{lorenz2014inertial,MoudafiOliny}) solve the more general maximal monotone inclusion problem:
\begin{eqnarray}
\text{Find } x\,\,\, \text{ s.t. }\, 0\in A(x)+B(x)\label{prob:mono}
\end{eqnarray}
where $A$ and $B$ are maximal monotone and $B$ is cocoercive\footnote{Setting $A=\partial g$ and $B=\nabla f$ recovers Prob.~(\ref{prob:1}).}. Other special cases were introduced as inertial versions of the KM iterations for finding fixed points \cite{mainge2008convergence,boct2015inertial}. In this paper we focus on convex optimization, which allows us to obtain less stringent convergence criteria than in those previous studies because we can use properties unique to convex functions. 
We note that GIPSA was originally suggested in \cite{liang2015activity}, however no theoretical properties were proven. GIPSA is also related (via discretization) to the continuous ODEs studied in \cite{attouch2015_2,su2014differential}.

I-FBS (i.e. the $\alpha_k=\beta_k$ case) includes as a particular parameter choice the Fast Iterative Shrinkage-Thresholding Algorithm (FISTA) introduced in \cite{Beck:2009} and inspired by earlier accelerated methods due to Nesterov \cite{nesterov1983method,nesterov2004introductory}, and G{\"u}ler \cite{guler1992new}. FISTA corresponds to I-FBS with the inertia sequence set to follow a specific rule which allows for a fast objective function convergence rate for Prob.~(\ref{prob:1}). Note that the recent variant of FISTA analyzed by Chambolle and Dossal \cite{chambolle2015convergence} and others \cite{attouch2015_2,attouch2015rate} is the only variant with a provably convergent iterate sequence. (Hereafter we refer to this variant as FISTA-CD due to the important contribution of Chambolle and Dossal).
Our global analysis when specialized to I-FBS requires $\alpha_k\leq\overline{\alpha}<1$, and therefore does not apply to any variant of FISTA. However our local analysis does apply to FISTA-CD.

The well-known property of FISTA and its variants is the ``fast" $O(1/k^2)$ objective function convergence rate for Prob.~(\ref{prob:1}). It is important to note that we do not expect this global objective function behavior to hold in general for I-FBS or GIPSA\footnote{In fact we expect the objective function values of GIPSA to reach the minimum with speed $o(1/k)$ if the parameters satisfy the conditions of Theorem \ref{thm:fista_lyap}. However this analysis is beyond the scope of this paper.}. Nevertheless the goal of this paper is not to study objective function convergence rates, but convergence of the iterates $\{x^k\}_{k\in\mathbb{N}}$ themselves, which is also important in practice \cite[p.~5]{prox_signalProcessing}. When we do compute convergence rates in the local analysis, they are asymptotically linear rates applicable to the iterates, i.e. $\|x^k-x^*\|\leq Cq^k$ for sufficiently large $k$, where $x^*$ is an optimal solution, and $q\in(0,1)$. One of the main findings of our local analysis for Prob.~(\ref{prob:sparse}) is that despite the optimal global sublinear convergence rate of FISTA, its local convergence performance can be greatly improved. This is important for applications where a high accuracy solution is needed, such as medical imaging \cite{choi2010compressed,lustig2007sparse}. 

\subsubsection{Global Analysis}
In our global analysis we aim to establish conditions on $\{\alpha_k,\beta_k,\lambda_k\}_{k\in\mathbb{N}}$ that imply the global weak convergence of the iterates $\{x^k,y^k,z^k\}_{k\in\mathbb{N}}$ of GIPSA to a solution of Prob.~(\ref{prob:1}). To the best of our knowledge no theoretical convergence study of (\ref{IFBS_yupdate})--(\ref{IFBS_xupdate}) in its full generality exists.
Special cases of GIPSA corresponding to different parameter choices have been studied previously in \cite{lorenz2014inertial,MoudafiOliny,mainge2008convergence}. However these analyses were not specialized to Prob.~(\ref{prob:1}) and therefore impose stricter conditions on the stepsize and inertia parameter than developed here. 

Our global analysis builds on the investigation of the inertial proximal algorithm of \cite{alvarez2000minimizing}. This algorithm corresponds to GIPSA when the smooth function $f$ is not present. Essentially our global analysis extends \cite[Theorem 3.1]{alvarez2000minimizing} to the composite case. We show that a multistep Lyapunov energy function is nonincreasing and this allows us to establish finiteness of the sum of the squared increments, i.e. $\sum_{k\in\mathbb{N}}\|x^k - x^{k-1}\|^2<\infty$. This condition is also needed for the local analysis. Weak convergence then follows via Opial techniques adapted from \cite{MoudafiOliny}. 

\subsubsection{Local Analysis}
The forward-backward nature of I-FBS (i.e. the $\alpha_k=\beta_k$ case) renders it amenable to a local analysis for Problems (\ref{prob:sparse})--(\ref{prob:l1LS}). It has been observed that FBS obtains \emph{local linear convergence} for Prob.~(\ref{prob:sparse}) and others \cite{Hale:2008,bredies2008linear,liang2014local,agarwal2012,liang2014convergence}. 
By this it is meant that after finitely many iterations, the iterates are permanently confined to a manifold containing the solution with respect to which the objective function is smooth. Thus after a finite period, convergence to a solution is linear if the local, smooth part of the function is also strongly convex, or a strict complementarity condition holds \cite{liang2014local,Hale:2008}. For Prob.~(\ref{prob:sparse}) the objective function is smooth with respect to vectors of fixed sign and support. 

We extend these results to I-FBS (including FISTA-CD). We show that I-FBS achieves local linear convergence and we determine the convergence rate in terms of the local curvature, the stepsize, and the inertia parameter. 
Importantly our analysis shows that adding the inertia term allows for a far better asymptotic convergence rate than is achievable with FBS (or FISTA). The local analysis borrows from the framework developed in \cite{Hale:2008}, however extensive differences emerge in order to incorporate the inertia term. 

We note that our local analysis results and techniques differ from what was presented in \cite{tao2015local}, which used a spectral analysis to study the local behavior of FBS and FISTA applied to Prob.~(\ref{prob:l1LS}). In contrast our analysis is based around exploiting the contractive properties of the \emph{soft-thresholding operator}, which is the proximal operator with respect to the $\ell_1$-norm. The authors of \cite{tao2015local} claim that both algorithms obtain local linear convergence when the minimizer is unique and a strict complementarity condition holds. 
Some of our results require neither of these conditions (Thms. \ref{thm:finite} and \ref{cor:q_lin}) while others depend on either strict complementarity (Thm. \ref{thm:furtherconv}) or solution uniqueness (Cor.~\ref{cor:upperb}). 
Unlike \cite{tao2015local}, we can compute Q-linear and R-linear convergence rates and this allows us to determine the optimal value for the inertia parameter. Many of our results also hold for the more general Prob.~(\ref{prob:sparse}). Our local analysis is also related to \cite{hare2004identifying} and we discuss this relationship in more detail in Sec.~\ref{sec:fistal1ls_finite}. 

We note that it is possible to derive upper bounds on the number of iterations not confined to the optimal smooth manifold within our analysis framework. To the best of our knowledge this is not possible in the competing frameworks \cite{tao2015local,hare2004identifying}. In some situations these upper bounds might be useful, however  in general they appear to be overly pessimistic compared to what is observed in practice. For now we do not attempt to provide tighter upper bounds. 



The FISTA variant FISTA-CD has stronger properties than the original version of \cite{Beck:2009} (in particular, see (\ref{abiggy})). In \cite{chambolle2015convergence,attouch2015_2} these properties were used to establish convergence of the iterates of FISTA-CD as well as the fast $O(1/k^2)$ objective function convergence rate. Recently in  \cite{attouch2015rate} Attouch and Peypouquet improved this to $o(1/k^2)$. In this paper we use these strong properties to establish the local convergence behavior of FISTA-CD for Problems (\ref{prob:sparse}) and (\ref{prob:l1LS}). We prove that FISTA-CD, exactly like I-FBS, obtains finite manifold identification for these problems. 
Furthermore, we show that after finitely many iterations FISTA-CD reduces to the form of a linear iterative system that has been studied previously in \cite{o2012adaptive}, allowing us to determine the asymptotic linear convergence rate. This rate is worse than that of the best choice for the inertia parameter in I-FBS and is comparable with the rate of (non-inertial) FBS. We then show that an adaptive restart scheme can be incorporated into FISTA-CD to obtain the optimal\footnote{among first-order methods} asymptotic convergence rate. Furthermore, unlike the optimal fixed choice of the inertia parameter, the restart scheme does not require knowledge of the local curvature parameter. Our restart modification of FISTA-CD preserves the optimal global convergence rate of FISTA-CD while also obtaining the optimal local convergence rate.


We note that restart techniques have been proposed before for accelerated methods, as well as conjugate gradient schemes, in the context of smooth and strongly convex problems \cite{o2012adaptive,su2014differential,monteiro2012adaptive}, \cite[p.~140]{bertsekas1999nonlinear}. It has been conjectured that restarting could improve the performance of FISTA even in the presence of nonsmooth regularizers \cite[\S 5.2]{o2012adaptive},\cite[p. 36]{cevher2014convex}. Our contribution is to show that this is indeed true for the case of the $\ell_1$-norm and to derive explicit convergence rates. 



The rest of the paper is organized as follows. In Section \ref{sec:prelim}, notation, definitions, assumptions and some preparatory results are presented. In Sections \ref{sec:fistalyapanal} through \ref{sec_discuss_param} we detail the conclusions of our global analysis of GIPSA for Prob.~(\ref{prob:1}). In Sections \ref{sec:l1ls} through \ref{sec:asFISTA} we give the results of our local convergence analysis of I-FBS and FISTA-CD for Problems (\ref{prob:sparse})--(\ref{prob:l1LS}). In Sec.~\ref{sec:sims}, a small synthetic numerical experiment on Prob.~(\ref{prob:l1LS}) is presented in order to corroborate some of our theoretical findings. Finally the proofs of all our results are given in Sections \ref{proofOfLemConstGrad1} through \ref{proofOfTheorem55}.

\section{Preliminaries}
\label{sec:prelim}
\subsection{Notation} 
For the most part the notation and conventions follow \cite{bauschke2011convex}. Thus $\mathcal{H}$ is a Hilbert space over the reals, $\langle\cdot,\cdot\rangle$ is the inner product and $\|\cdot\|$ is the induced norm. For $\mathbb{R}^n$ we assume the standard Euclidean norm and inner product and use $\|\cdot\|_1$ to denote the $\ell_1$-norm. The notation $\mathbb{R}_+$ denotes the set of all nonegative reals. Let $\Gamma_0(\mathcal{H})$ be the set of all closed, convex and proper functions from $\mathcal{H}$ to $(0,\infty]$. For any $g:\mathcal{H}\to(0,\infty]$ and point $x\in\mathcal{H}$, we denote by $\partial g(x)$ the \emph{subdifferential} at $x$ \cite[Def. 16.1]{bauschke2011convex}, and by $\text{dom}\,\,\partial g\subset \mathcal{H}$ the set of $x$ such that $\partial g(x)$ is nonempty.

For $a:\mathbb{R}\to\mathbb{R}$, $b:\mathbb{R}\to\mathbb{R}$, and $c\in[-\infty,+\infty]$, the notation $a(l)=O(b(l))$ (resp. $a(l)=\Omega(b(l))$)  means there exists a constant $C\geq0$ such that $\limsup_{l\to c}|a(l)/b(l)|\leq C$ (resp. $\liminf_{l\to c}|a(l)/b(l)|\geq C$). 
We will say a sequence $\{x^k\}_{k\in\mathbb{N}}\subset\mathcal{H}$ converges \emph{R-linearly} to $x^*\in\mathcal{H}$ with rate of convergence $q\in(0,1)$, if $\|x^k-x^*\|=O(q^k)$.  We say $x^k$ converges to $x^*$ \emph{Q-linearly} with rate $q\in(0,1)$ if
$
\lim_{k\to\infty}\left\{\|x^k - x^*\|/\|x^{k-1}-x^*\|\right\} = q
$.
Colloquially we will refer to both Q-linear and R-linear convergence simply as linear convergence. Occasionally we say local or asymptotically linear convergence. 
We use $x^k\to x^*$ to denote strong convergence and $x^k\rightharpoonup x^*$ to denote weak convergence. 

For Prob.~(\ref{prob:1}) define the \emph{optimal value} as
$
F^* \triangleq \inf_{x\in\mathcal{H}} F(x)
$
and the \emph{solution set} as
$
X^* \triangleq \{x\in\mathcal{H}:F(x) = F^*\}
$
which may be empty.
For the sequence $\{x^k\}_{k\in\mathbb{N}}$ generated by (\ref{IFBS_yupdate})--(\ref{IFBS_xupdate}), let $\Delta_{k}$ denote $x^{k}-x^{k-1}$ for all $k\in\mathbb{N}$. Given a function $a:\mathbb{R}\to\mathbb{R}$, we say that the \emph{iteration complexity} of a method for minimizing $F$ is $\Omega\lbr a(\epsilon)\rbr$ if $k=\Omega\lbr a\lbr\epsilon\rbr\rbr$ implies $F(x^k)-F^*=O(\epsilon)$ as $\epsilon\to 0$. 

For a sorted set $S\subseteq\{1,2,\ldots,n\}$ with no repeated elements, let $S(i),i=1,\ldots,|S|$ be the $i$th element of $S$. For a matrix $A\in\mathbb{R}^{m\times n}$, $A_S$ will denote the matrix in $\mathbb{R}^{m\times |S|}$ formed by taking the columns corresponding to the elements of $S$. That is $A_S(i,j)=A(i,S(j))$. For a vector $v\in\mathbb{R}^{n}$, $v_S$ will denote the $|S|\times 1$ vector with entries given by $v_S(i)=v(S(i))$. The notation $(v_S,0)$ will denote the vector in $\mathbb{R}^n$ whose $j$th entry is $v(j)$ if $j\in S$ and $0$ otherwise. 
The range space and null space of a matrix $A$ are denoted by $\mathcal{R}(A)$ and $\mathcal{N}(A)$ respectively.
Given $c\in\mathbb{R}$ and $x\in\mathbb{R}^n$, $\sgn(c)$ is defined as $+1$ if $c\geq0$ and $-1$ if $c<0$, $\sgn(x)$ is simply applying $\sgn(\cdot)$ elementwise. Finally  $[c]_+\triangleq\max(c,0)$. 

The following identity appears in many convergence analyses and we will use it many times. For all $x,y,z\in\mathcal{H}$, 
\begin{eqnarray} 
\langle x-y,x-z\rangle = \frac{1}{2}\|x-y\|^2+\frac{1}{2}\|x-z\|^2-\frac{1}{2}\|y-z\|^2\label{biglem}.
\end{eqnarray}

\subsection{Properties of Convex and Smooth Functions}
Now we list some properties of the subdifferential, as well as convex and smooth functions. 
For the Fr\'echet and  G\^ateaux definitions of differentiability we refer to \cite[Definition 2.45 and 2.43]{bauschke2011convex}. Note that Fr\'echet differentiability on a neighborhood of a point implies G\^ateaux differentiability at that point, and the two derivatives agree \cite[Lemma 2.49(i)]{bauschke2011convex}.  
Consider $f:\mathcal{H}\to(-\infty,+\infty]$. Then 
\begin{eqnarray}
\label{eq:ineq2}
\langle t,u-v\rangle &\geq& f(u)-f(v),\quad\forall v\in\mathcal{H},\,\,u\in\text{dom}\,\partial f,\,\text{ and }\,t\in\partial f(u),
\end{eqnarray}
and
\begin{eqnarray} 
\langle t-p,u-v\rangle \geq 0\label{montone_op},\quad\forall u,v\in\text{dom}\,\partial f,\,\, t\in\partial f(u) \text{ and } p\in\partial f(v).
\end{eqnarray}
For a proper and convex function which is G\^ateaux differentiable everywhere on $\mathcal{H}$, (\ref{eq:ineq2})--(\ref{montone_op}) hold for all $u,v\in\mathcal{H}$ \cite[Prop.~17.10]{bauschke2011convex} and $\partial f(x)=\{\nabla f(x)\}$  (i.e. a singleton) everywhere. 

We say that a Fr\'echet differentiable function $f$ has $L$-Lipschitz continuous gradient if
$
\|\nabla f(y)-\nabla f(x)\|\leq L\|y-x\|,\ \forall x,y\in\mathcal{H}.
\label{eq:lips}
$
For such a function \cite[Thm. 18.15 (iii)]{bauschke2011convex}: 
\begin{eqnarray}
\label{eq:ineq1}
f(u)-f(v)&\leq&\langle\nabla f(v),u-v\rangle + \frac{L}{2}\|u-v\|^2,\quad\forall u,v\in\mathcal{H}.
\end{eqnarray} 
The gradient $\nabla f$ of a convex and Fr\'echet differentiable function is $L$-Lipschitz continuous if and only if \cite[Cor.~18.16]{bauschke2011convex}
\begin{eqnarray}
\langle\nabla f(u)-\nabla f(v),u-v\rangle 
\geq 
\frac{1}{L}\|\nabla f(u)-\nabla f(v)\|^2
,\quad\forall u,v\in\mathcal{H}.
\label{def_coercive}
\end{eqnarray}
This is the celebrated Baillon-Haddad Theorem. 

\subsection{Proximal Operators}
\label{sec:prox}
The proximal operator $\prox_g:\mathcal{H}\to\mathcal{H}$ with respect to a function $g\in\Gamma_0(\mathcal{H})$ is defined implicitly as:
$
y-\prox_g(y) \in \partial g(\prox_g(y)),
\label{implicit}
$
and explicitly as
$$
\prox_g(y)=\arg\min_{x\in\mathcal{H}}\left\{\frac{1}{2}\|x-y\|^2+ g(x)\right\},\quad\forall y\in\mathcal{H}. 
\label{eq:prox_def}
$$
The proximal operator is a well-defined mapping from $\mathcal{H}$ to $\text{dom}\,\partial g$ \cite[Prop. 23.2, Example 23.3]{bauschke2011convex}. In light of the implicit definition of the proximal operator we point out that the update equation for GIPSA given in (\ref{IFBS_xupdate}) can be written implicitly as 
\begin{eqnarray}
0\in x^{k+1} - y^{k+1} +\lambda_k\partial g(x^{k+1}) + \lambda_k\nabla f(z^{k+1}).\label{implicitupdate}
\end{eqnarray}
Now $\rho\|\cdot\|_1\in\Gamma_0(\mathbb{R}^n)$ and the proximal operator associated with it is the shrinkage and soft-thresholding operator $S_{\rho}(v):\mathbb{R}\to\mathbb{R}$, applied elementwise. It is defined as
\begin{eqnarray}
\label{def_softThresh}
S_{\rho}(v)\triangleq\left[|v|-\rho\right]_+\sgn(v)
\,\,\,\text{and}\,\,\,
\label{soft_thresh_def}
\{\prox_{\rho\|\cdot\|_1}(z)\}_i=S_{\rho}(z_i),\, i=1,2,\ldots,n.
\end{eqnarray} 

\subsection{Assumptions and Optimality Conditions}
Now we are ready to precisely state the assumptions used throughout the paper.
\label{sec_assump_optimal}
\vspace{0.1cm}
\begin{assump}[Problems (\ref{prob:1})--(\ref{prob:sparse})]
	The functions $f$ and $g$ are in $\Gamma_0(\mathcal{H})$, $\text{dom}\,\partial g$ is nonempty, $f$ is Fr\'echet differentiable everywhere and has an $L$-Lipschitz continuous gradient with $L>0$, and $F^*>-\infty$.
\end{assump}
\vspace{0.1cm}

The optimality conditions for Prob.~(\ref{prob:1}) under Assumption 1 are as follows. A vector $x^*\in X^*$ if and only if \cite[Corollary 26.3 (vi)]{bauschke2011convex}
\begin{eqnarray}
0\in\partial F(x^*)=(\partial g +\nabla f)(x^*) = \partial g(x^*) +\{\nabla f(x^*)\}.
\label{optcond}
\end{eqnarray}
Note that this is equivalent to $x^*$ satisfying
\begin{eqnarray}
x^* = \prox_{\lambda g}(x^* - \lambda\nabla f(x^*))
\label{optcond_fixed}
\end{eqnarray}
for all $\lambda>0$ \cite[Corollary 26.3 (viii)]{bauschke2011convex}. Thus $x^*$ is a solution to Prob.~(\ref{prob:1}) if and only if it is a fixed point of the \emph{forward-backward operator}: 
$
T_\lambda(x)\triangleq\prox_{\lambda g}(x-\lambda\nabla f(x)).
$
Note that $T_{\lambda}$ is nonexpansive so long as $0\leq\lambda< 2/L$ \cite[Thm. 25.8]{bauschke2011convex}


%

\label{sec:propsparse}




The function $\frac{1}{2}\|Ax-b\|^2$ is differentiable and has gradient equal to $A^\top (Ax-b)$ which is Lipschitz continuous with Lipschitz constant equal to the largest eigenvalue of $A^\top A$. The objective function in Prob.~(\ref{prob:l1LS}) is bounded below by $0$. As previously stated, $\rho\|\cdot\|_1\in\Gamma_0(\mathcal{H})$ and $\text{dom}\,\,\partial\|\cdot\|_1=\mathbb{R}^n$. Therefore Prob.~(\ref{prob:l1LS}) satisfies Assumption 1. Thus results proved for Prob.~(\ref{prob:1}) hold for all problems, while results proved for Prob.~(\ref{prob:sparse}) also hold for Prob.~(\ref{prob:l1LS}). 
Note that the solution set $X^*$ of Prob.~(\ref{prob:l1LS}) is always nonempty.

\subsection{Properties of the Solution Set of Prob.~(\ref{prob:sparse})}\label{props_of_l1}
\begin{lemma}
	Suppose Assumption 1 holds for Prob.~(\ref{prob:sparse}) and $X^*$ is nonempty, then there exists a vector $h^*\in\mathbb{R}^n$ such that for all $x^*\in X^*$, $\nabla f(x^*)=h^*$. Furthermore, for all $i\in\{1,2,\ldots,n\}$ and $x^*\in X^*$: $-h_i^* x_i^*\geq 0$. Finally
	\begin{eqnarray*}
		\frac{h_i^*}{\rho}\left\{
		\begin{array}{ll}
			=-1:&\hbox{ if }\, \exists\ x^*\in X^*:x^*_i>0 \\
			=+1:& \hbox{ if }\, \exists\ x^*\in X^*:x^*_i<0\\
			\in [-1,1]:&\quad \text{else}.
		\end{array}
		\right.
	\end{eqnarray*}
	\label{thm:constGrad}
\end{lemma}

\emph{Proof.} See Sec.~\ref{proofOfLemConstGrad}.

Let
$
\label{Edef}
E\triangleq\{i:|h_i^*|=\rho\}
$
and note that
$
\label{Ddef}
E^c = \{i:|h_i^*|<\rho\}.
$
Throughout the paper we will assume the elements of $E$ are in increasing order. By Lemma \ref{thm:constGrad}, we infer that $\supp(x^*)\subseteq E$ for all $x^*\in X^*$. The set $E$ will be crucial to our local analysis. 



\subsection{Properties of FISTA-CD}
Chambolle and Dossal \cite{chambolle2015convergence} analyzed a variant parameter choice of FISTA which has the $O(1/k^2)$ global objective function convergence rate and also convergence of the sequence $\{x^k\}_{k\in\mathbb{N}}$ to a minimizer (see also \cite{attouch2015_2,attouch2015rate}). They considered the following parameter choice for GIPSA (more specifically I-FBS), which we refer to as FISTA-CD:
\begin{eqnarray}
\label{chamchoice}
x^1=x^0,\quad \lambda_k=\lambda\in(0,1/L],\quad \alpha_{k} = \beta_k=\frac{k-1}{k+a},\quad a>2,\quad\forall k\in\mathbb{N}. 
\end{eqnarray} 
For a discussion on how to choose $a$ see \cite[\S 4]{chambolle2015convergence}. The choice $a=2$, which is not permitted in the analysis of \cite{chambolle2015convergence}, corresponds to one of the classical parameter choices for Nesterov's accelerated method for which convergence of the iterates is still unknown \cite[Remark 2.17]{attouch2015_2}. We now detail the important properties of FISTA-CD derived in \cite{chambolle2015convergence} which we need for our analysis.
\begin{lemma}[\cite{chambolle2015convergence}]
	\label{lemCham}
	Suppose Assumption 1. holds for Prob.~(\ref{prob:1}), $X^*$ is nonempty, and $\{\lambda_k\}_{k\in\mathbb{N}}$ and $\{\alpha_k\}_{k\in\mathbb{N}}$ are chosen as in (\ref{chamchoice}). Then for the iterates $\{x^k\}_{k\in\mathbb{N}}$ of (\ref{IFBS_yupdate})--(\ref{IFBS_xupdate}): 
	\begin{enumerate}
		\item \cite[Theorem 4.1: Eq. (25)]{chambolle2015convergence} 
		\begin{eqnarray}
		\sum_{k=1}^\infty\sum_{j=1}^k \left(\prod_{l=j}^k\alpha_l\right)\|x^j-x^{j-1}\|^2 <\infty.\label{abiggy}
		\end{eqnarray}
		\item \cite[Theorem 4.1]{chambolle2015convergence} There exists $\hat{x}\in X^*$ such that $x^k\rightharpoonup \hat{x}$.  
	\end{enumerate}
\end{lemma}

\section{Main Results}
\label{sec_mainResults}
\subsection{Global Convergence Analysis of GIPSA}
\label{sec:fistalyapanal}
In this section we state conditions on $\{\alpha_k,\beta_k,\lambda_k\}_{k\in\mathbb{N}}$ which imply weak convergence of the iterates $\{x^k,y^k,z^k\}_{k\in\mathbb{N}}$ of (\ref{IFBS_yupdate})--(\ref{IFBS_xupdate}) to a minimizer of Prob.~(\ref{prob:1}) under Assumption 1. These conditions also imply finite summability of the squared increments of the sequence, which will be useful in the local analysis. The finite summability result also makes it trivial to prove criticality of the limit points which we include for completeness.
\label{sec:lyap}

\begin{theorem}\label{thm:fista_lyap}
	For Prob.~(\ref{prob:1}), suppose Assumption 1. holds. Assume $\{\lambda_k\}_{k\in\mathbb{N}}$ is positive and nondecreasing, and there exists $\epsilon>0$, $0<\gamma<2$ and $0\leq\overline{\beta}<1$ such that sequences $ \{\lambda_k,\alpha_k,\beta_k\}_{k\in\mathbb{N}}$ satisfy:   
	\begin{eqnarray}
	&& 0\leq\alpha_k\leq 1, 
	\quad
	0\leq\beta_k\leq\overline{\beta},
	\quad
	\lambda_k\alpha_k \leq \frac{\beta_k}{L},
	\quad
	\lambda_k\leq \frac{2-\gamma}{L}
	\nonumber\\&&\label{param_conditions}
	\,\,\quad\text{ and }\quad
	2- \lambda_k L (1-\alpha_k)-\beta_k -\beta_{k+1}\geq \epsilon
	\end{eqnarray}
	for all $k\in\mathbb{N}$. Then for the iterates $\{x^k,y^k,z^k\}_{k\in\mathbb{N}}$ of (\ref{IFBS_yupdate})--(\ref{IFBS_xupdate}):
	\begin{enumerate}[(i)]
		\item $\sum_{k\in\mathbb{N}}\|x^k - x^{k-1}\|^2 < \infty$,\label{state:finitesum}
		\item $d(0,\partial F(x^k))\to 0$ as $k\to\infty$.\label{criticality}
		\item If $X^*$ is nonempty then there exists $\hat{x}\in X^*$ such that $x^k\rightharpoonup \hat{x}$, $y^k\rightharpoonup \hat{x}$ and $z^k\rightharpoonup \hat{x}$. \label{state:weak}
	\end{enumerate}
\end{theorem}

\emph{Proof.}
See Sec.~\ref{proofOfThm41}.

With some effort Theorem \ref{thm:fista_lyap} can be extended to inexact proximal operators through the use of the $\epsilon$-enlarged subdifferential under a summability condition on the errors \cite{burachik1999varepsilon}. It can also be extended to versions which incorporate a relaxation parameter. 
To simplify the presentation, proof, and notation, we do not detail these elaborations. 

For the special case where $\alpha_k=0$, Theorem \ref{thm:fista_lyap} provides more general parameter constraints than existing guarantees  derived in \cite{MoudafiOliny}. Suppose $\lambda_k=\lambda\in[0,2/L)$, then \cite{MoudafiOliny} requires $\beta_k$ to be nondecreasing and to satisfy $0\leq\beta_k\leq\overline{\beta}$ where 
$
\overline{\beta}< (2-\lambda L)/6.
$
On the other hand, Theorem \ref{thm:fista_lyap} requires: 
$
\beta_k + \beta_{k+1}\leq 2-\lambda L -\epsilon,
$
which is satisfied if $\overline{\beta}< (2-\lambda L)/2$. Note that \cite{MoudafiOliny} and Theorem \ref{thm:fista_lyap} have the same requirement on the stepsize.

\subsection{Specialized Conditions for I-FBS}\label{spec_IFBS}
We now simplify the conditions for the case of I-FBS, i.e. $\alpha_k=\beta_k$. For consistency, let $\overline{\alpha}=\overline{\beta}$. In this case:
$$
2-\alpha_k - \alpha_{k+1} + \lambda_k L(\alpha_k-1)
\geq 
1-\alpha_{k+1}\geq 1-\overline{\alpha}>0.
$$
Therefore $\epsilon=1-\overline{\alpha}$ satisfies (\ref{param_conditions}). Next note that if we choose any $\gamma<1$, the condition on the stepsize simplifies to $\lambda_k\leq 1/L$ for all $k$. In the case of FBS: i.e. $\alpha_k=\beta_k=0$, Thm. \ref{thm:fista_lyap} allows for larger stepsizes: $\lambda_k L\leq 2-\gamma < 2$, which agrees with the standard criteria for FBS (e.g. \cite[Thm. 25.8]{bauschke2011convex}). We formalize this in the following corollary.
\begin{corollary}
	Assume $\{\lambda_k\}_{k\in\mathbb{N}}$ is nondecreasing, $\lambda_k\in(0,1/L]$, $\alpha_k=\beta_k$, and $0\leq\alpha_k\leq\overline{\alpha}<1$ for all $k$. Then for the iterates of GIPSA for Prob.~(\ref{prob:1}), Assumption 1 implies (\ref{state:finitesum}) and (\ref{criticality}) of Theorem \ref{thm:fista_lyap}. Assumption 1 and nonemptiness of $X^*$ imply (\ref{state:weak}) of Theorem \ref{thm:fista_lyap}.\label{cor_IFBS_weak}
\end{corollary}

Note that the condition on $\alpha_k$ is more general than the requirement on the inertia parameter given in \cite{lorenz2014inertial} which is 
$$
1-3\alpha_k - \lambda L(1-\alpha_k)^2/2 \geq \eta
$$
where $\lambda_k=\lambda\in(0,2/L]$ for all $k$ and $\eta>0$ is some constant. Note that \cite{lorenz2014inertial} does allow larger values of the stepsize $\lambda$, up to $2/L$ so long as $\alpha_k$ is sufficiently small. 

We emphasize that Corollary \ref{cor_IFBS_weak} does not apply to any of the FISTA variants because in all such algorithms $\alpha_k\to 1$. See \cite{chambolle2015convergence} for a proof of weak convergence of the iterates of FISTA-CD. 

It is interesting to note that for FBS the convergence criteria are the same for Prob.~(\ref{prob:mono}) (monotone inclusion problem) and Prob.~(\ref{prob:1}) \cite[Thm. 25.8 and Cor.~27.9]{bauschke2011convex}. However for GIPSA and I-FBS, this does not appear to be the case.


\subsection{Discussion of the General Case}\label{sec_discuss_param}
We have discussed the special cases: $\alpha_k=\beta_k$ and $\alpha_k=0$, we now discuss the general case. To simplify the discussion, consider fixed choices, i.e. $\{\alpha_k,\beta_k,\lambda_k\}=\{\alpha,\beta,\lambda\}$ for all $k$. 
Then (\ref{param_conditions}) becomes
\begin{eqnarray}\label{fixed_param_conds}
\alpha\in[0,1],\quad\beta\in[0,1),
\quad
0<
\lambda L
&\leq&
\min\left\{\frac{\beta}{\alpha},\frac{2(1-\beta)-\epsilon}{1-\alpha}\right\}
\end{eqnarray}
for some $\epsilon>0$ with the convention: $0/0=\infty$. Now if we set $\epsilon$ to $0$, the two arguments to $\min$ in (\ref{fixed_param_conds}) are equal if
$
\alpha=\alpha^*(\beta) = \frac{\beta}{2-\beta}.
$
Substituting this into the expression yields $\lambda L < 2-\beta$. If $\alpha<\alpha^*(\beta)$ then the right-hand expression in the argument of $\min$ is the smallest, else it is the left-hand expression. Thus the condition on $\lambda$ is 
\begin{eqnarray*}
	\lambda L < 
	\left\{
	\begin{array}{cc}
		\frac{2(1-\beta)}{1-\alpha}: & \text{if }0\leq\alpha\leq \alpha^*(\beta)\\
		\frac{\beta}{\alpha}: & \text{if }\alpha^*(\beta)\leq\alpha\leq 1.
	\end{array}
	\right.
\end{eqnarray*}

While $\alpha=\alpha^*(\beta)$ provides the largest range of feasible stepsize parameters according to our theoretical convergence analysis, we do not claim that it is the ``best" choice for a given instance of Prob.~(\ref{prob:1}). For I-FBS for Prob.~(\ref{prob:sparse}) our local convergence analysis derives some ``good" parameter choices (See Sections \ref{sec:l1ls}--\ref{sec:asFISTA}). Some of these choices have both good local and global convergence properties. However determining good parameter choices more generally for GIPSA is beyond the scope of this paper. Nevertheless it is important to establish general conditions for convergence before attempting to determine appropriate choices via an empirical study or theoretical analysis.


\subsection{Finite Convergence Results for I-FBS}
\label{sec:l1ls}\label{sec:convFistal1ls}
We now turn our attention to Problems (\ref{prob:sparse})--(\ref{prob:l1LS}) and establish the local convergence behavior of I-FBS and FISTA-CD. 
\label{sec:fistal1ls_finite}
The upcoming theorem proves convergence in a finite number of iterations for the components in $E^c$ to $0$, and for the components in $E$ to the optimal sign (recall $E\triangleq\{i:|h_i^*|=\rho\}$ where $h^*$ is defined in Lemma \ref{thm:constGrad}). Following the terminology of \cite{hare2004identifying,liang2014local} we will refer to this as the ``finite active manifold identification" property. The manifold in the $\ell_1$-norm setting is the halfspace of vectors with support a subset of $E$ and nonzero components with sign equal to $-h_i^*/\rho$. 

\begin{theorem}
	For Prob.~(\ref{prob:sparse}) suppose that Assumption 1 holds and $X^*$ is nonempty, thus there exists $h^*\in\mathbb{R}^n$ satisfying the conditions of Lemma \ref{thm:constGrad}. Assume that either:
	\begin{enumerate}
		\item $\{\lambda_k\}_{k\in\mathbb{N}}$ is nondecreasing, $0<\lambda_k\leq1/L$, $\alpha_k=\beta_k$ and $0\leq\alpha_k\leq \overline{\alpha}<1$ for all $k\in\mathbb{N}$, or
		\item $\{\alpha_k,\lambda_k\}_{k\in\mathbb{N}}$ are chosen according to (\ref{chamchoice}) (i.e.  FISTA-CD),
	\end{enumerate}
	then for all but finitely many $k$ the iterates $\{x^k,y^k\}_{k\in\mathbb{N}}$ of (\ref{IFBS_yupdate})--(\ref{IFBS_xupdate}) satisfy
	\begin{eqnarray}
	\label{thm:finite:eq2}
	\sgn\left(y_i^{k}-\lambda_{k-1}\nabla f(y^{k})_i\right)
	&=&-\frac{h^*_i}{\rho},\ \forall i:|h^*_i|=\rho,
	\end{eqnarray}
	and
	\begin{eqnarray}
	x^k_i =y^k_i= 0,\ \forall i:|h^*_i|<\rho.
	\label{thm:finite:eq1}
	\end{eqnarray}
	\label{thm:finite}
\end{theorem}

\emph{Proof.}
See Sec.~\ref{proofOfTheroem5}.

Note that if  $x_i^{k}\neq 0$, then (\ref{def_softThresh}) implies that $\sgn(x_i^{k})=\sgn(y_i^{k}-\lambda_k\nabla f(y^{k})_i)$. Note that \cite[Theorem 4.5]{Hale:2008} is recovered when $\alpha_k=0$ for all $k$. 


The authors of \cite{hare2004identifying} studied finite convergence results for prox-regular and partially smooth functions, which includes Prob.~(\ref{prob:sparse}). Specialized to this problem, the analysis of \cite[Theorem 5.3]{hare2004identifying} establishes finite convergence in support and sign for any algorithm which produces a convergent iterate sequence, under the following additional condition: $0\in \text{rint}(\partial F(x^*))$ for the limit $x^*=\lim_{k\to\infty}x^k$. In the context of Prob.~(\ref{prob:sparse}) this condition is equivalent to the ``strict complementarity condition" discussed in Sec.~\ref{sec:locallin2}, i.e. $E=\supp(x^*)$. 
In contrast, Theorem \ref{thm:finite} is more general in that it proves finite convergence to $0$ on $E^c\subseteq\supp(x^*)^c$ and sign on $E$. It does not require $E=\supp(x^*)$. However when this is true, Theorem \ref{thm:finite} coincides with \cite[Theorem 5.3]{hare2004identifying}. 

Given that I-FBS converges in a finite number of iterations to the optimal manifold, it could be desirable to switch to a local procedure which searches in the space of lower dimension. Indeed for Prob.~(\ref{prob:l1LS}), if the solution is unique and the support and sign of the solution are known, then the values of the nonzero entries can be computed by solving a linear system with dimension equal to the number of nonzero entries \cite[p.~20]{bach2011convex}. Theorem \ref{thm:finite} also motivates combining two-stage ``active-set" strategies such as the one described in \cite{wen2012convergence} with I-FBS or FISTA-CD. Active-set strategies alternate between iterated shrinkage-thresholding updates to identify the active manifold, and local optimization procedures to estimate the nonzero entries. Using I-FBS/FISTA-CD to identify the active manifold within such a framework is an interesting topic for future work. 

\subsection{Reduction to Smooth Minimization}
Theorem \ref{thm:finite} allows us to characterize the behavior of I-FBS after a manifold identification period of finite duration. In the following theorem, we show that after a finite number of iterations, I-FBS (including FISTA-CD) reduces indefinitely to minimizing a smooth function over $E$ subject to an orthant constraint. 

\begin{theorem}
	\label{cor:q_lin}
	For Prob.~(\ref{prob:sparse}) suppose that Assumption 1 holds and $X^*$ is nonempty, thus there exists $h^*\in\mathbb{R}^n$ satisfying the properties of Lemma \ref{thm:constGrad}. Recall $E\triangleq\{i:|h_i^*|=\rho\}$ and let $\phi:\mathbb{R}^{|E|}\to\mathbb{R}$ be defined as
	\begin{eqnarray}
	\phi(x_E)&\triangleq&-(h_E^*)^\top x_E+f\left((x_E,0)\right),
	\label{eq:phidef}
	\end{eqnarray}
	where $x\in\mathbb{R}^n$. Let the set $O_E\subset \mathbb{R}^{|E|}$ be defined as
	\begin{eqnarray}
	O_E&\triangleq&\{v\in\mathbb{R}^{|E|}:-\sgn(h_{E(j)}^*)\, v_j\geq 0,\ \forall j\in \{1,2,\ldots,|E|\}\}.
	\label{O_E_def}
	\end{eqnarray} 
	Assume that either
	\begin{enumerate}
		\item $\{\lambda_k\}_{k\in\mathbb{N}}$ is nondecreasing, $0<\lambda_k\leq1/L$, $\alpha_k=\beta_k$ and $0\leq\alpha_k\leq \overline{\alpha}<1$, for all $k$, or
		\item $\{\alpha_k,\lambda_k\}_{k\in\mathbb{N}}$ are chosen according to FISTA-CD in (\ref{chamchoice}),
	\end{enumerate} 
	then, for all but finitely many $k$, the iterates $\{x^k,y^k\}_{k\in\mathbb{N}}$ of (\ref{IFBS_yupdate})--(\ref{IFBS_xupdate}) satisfy 
	\begin{eqnarray}
	x_E^{k+1}&=& P_{O_E}\lbr y_E^{k+1}-\lambda_k\nabla\phi(y_E^{k+1})\rbr,
	\label{eq:proj_HB2}
	\end{eqnarray}
	and $F(x^k) = \phi(x_E^k)$, where $F(x)=f(x)+\rho\|x\|_1$ and $P_{O_E}$ is the orthogonal projector onto $O_E$. 
\end{theorem}

\emph{Proof.}
See Sec.~\ref{Proof_q_lin}. The result of \cite[Corollary 4.6]{Hale:2008} is recovered when $\alpha_k=0$ for all $k$.

\subsection{Local Linear Convergence Under Local Strong Convexity}
\label{sec:paramsfist}
The analysis of the previous two sections shows that, after a finite number of iterations, I-FBS reduces to minimizing the function $\phi$ subject to an orthant constraint. This function can be strongly convex even if $f$ does not have this property. If $\phi$ is strongly convex, then local linear convergence can be achieved, as we prove in the following Corollary. Note that strong (in fact strict) convexity of $\phi$ implies solution uniqueness for Prob.~(\ref{prob:sparse}).

\begin{corollary}
	For Prob.~(\ref{prob:sparse}) suppose that Assumption 1 holds and $\phi$ defined in (\ref{eq:phidef}) is strongly convex. Let $l_E$ be the strong convexity parameter of $\phi$. If $\lambda\in(0,1/L]$, $0<\mu\leq l_E$,
	\begin{eqnarray}
	\label{corupper}
	\lambda_k=\lambda\hbox{ and } \beta_k=\alpha_k=\frac{1-\sqrt{\mu\lambda}}{1+\sqrt{\mu\lambda}}\quad\forall k\in\mathbb{N},
	\end{eqnarray}
	then the iterates $\{x^k\}_{k\in\mathbb{N}}$ of (\ref{IFBS_yupdate})--(\ref{IFBS_xupdate}) converge to the unique solution $x^*$ of Prob.~(\ref{prob:sparse}) R-linearly and $F(x^k)$ converges to $F^*$ R-linearly 	where $F(x)=f(x) + \rho\|x\|_1$. Specifically
	\begin{eqnarray*}
		\|x^k - x^*\|^2 = 	O\label{Frate}
		\left( 
		\left(1-\sqrt{\mu\lambda}\right)^{k}
		\right)\label{iterate_rate}
		\text{ and }\,
		F\lbr x^k\rbr-F^*
		= 
		O
		\left( 
		\left(1-\sqrt{\mu\lambda}\right)^{k}
		\right).
	\end{eqnarray*}
	\label{cor:upperb}
\end{corollary}
\emph{Proof.}
	Recall the definition of $E\triangleq\{i:|h_i^*|=\rho\}$. We consider $k$ large enough that I-FBS has reduced to minimizing the $l_E$-strongly convex function $\phi$, i.e. (\ref{eq:proj_HB2}) holds, $x_{E^c}^k=0$, and $F(x^k)=\phi(x_E^k)$.  The result can now be seen by considering Nesterov's constant momentum scheme of \cite[p.~76]{nesterov2004introductory}, however the variable $\mu$ now represents a lower bound for the true strong convexity parameter of $\phi$. It can be verified that this does not change the result given in \cite[Thm 2.2.3]{nesterov2004introductory}. Furthermore we allow stepsizes other than $1/L$, which is discussed on \cite[p.~72]{nesterov2004introductory}. Finally the minimization is with respect to the orthant $O_E$ defined in (\ref{O_E_def}). This simple modification of Nesterov's scheme is discussed in \cite[Algorithm (2.2.17)]{nesterov2004introductory}. 
	\qed

Note this local linear convergence result does not depend on strict complementarity (i.e. $E=\supp(x^*)$) unlike the local analysis of FBS in \cite{liang2014local,bredies2008linear}.
Suppose $\mu=l_E$ and $\lambda=1/L$ then the convergence rate and iteration complexity are respectively
\begin{eqnarray}\label{bnn}
F(x^k)-F^* = O\left(\left(1-\sqrt{\frac{l_E}{L}}\right)^k\right),\,\,
\text{iter. comp.}= \Omega\lbr\sqrt{\frac{L}{l_E}}\log\frac{1}{\epsilon}\rbr. 
\label{eq:itercomp}
\end{eqnarray}
Given the nature of $\phi$ this iteration complexity is optimal \cite{nesterov2004introductory}. Indeed it is better than the iteration complexity of FBS \cite{Hale:2008} (which corresponds to I-FBS with $\alpha_k$ equal to $0$) which is
$
\label{istaIterComp}
\Omega\lbr(L/l_E)\log1/\epsilon\rbr.
$



Other parameter choices, such as Constant Scheme III of \cite[p.~84]{nesterov2004introductory}, will also achieve local linear convergence with the same rate. 
However these choices along with (\ref{corupper}) are difficult to use in practice as they depend on $l_E$, which is hard to estimate. In Sec.~\ref{sec:asFISTA} we will show how the rate and corresponding iteration complexity in (\ref{eq:itercomp}) can be achieved without knowledge of $l_E$ by combining a restart scheme with FISTA-CD.

  \subsection{Local Linear Convergence Under Strict Complementarity}
  \label{sec:locallin2}
  Local linear convergence can also be proved for Prob.~(\ref{prob:l1LS}) without requiring solution uniqueness. We require $\lim_{k\to\infty} x^k \triangleq x^*\in X^*$ to obey the so-called ``strict complementarity" condition: 
  $
  E=\supp(x^*),
  $
  where $E\triangleq\{i:|h_i^*|=\rho\}$. This is a common assumption also used in \cite{Hale:2008,liang2014local,tao2015local,hare2004identifying,bredies2008linear}. Note that this condition is not necessary for $x^*$ to be the unique minimizer for Prob.~(\ref{prob:sparse}) \cite[Example (4)]{zhang2015necessary}. It is also not sufficient, which can be seen by considering the following instance of Prob.~(\ref{prob:l1LS}) taken from  \cite[Example (4)]{zhang2015necessary}:
  \begin{eqnarray*} 
  A
  =
  \left[
  \begin{array}{ccccccc}
  1 & & & 0 & & &2\\
  0 & & & 2 & & &-2
  \end{array}
  \right]
  ,
  b=
  \left[
  \begin{array}{c}
  1.5\\
  1
  \end{array}
  \right],
  \rho=1.
  \end{eqnarray*}
  This example has $E=\{1,2,3\}$ and has infinitely many solutions which satisfy strict complementarity, such as $(1/4,3/8,1/8)^\top$. The name ``strict complementarity" comes from considering the dual problem to Prob.~(\ref{prob:sparse}) \cite[\S 6]{nesterov2013gradient}.
  
  
  First we state the following proposition which shows that the proximal step (\ref{IFBS_xupdate}) of I-FBS reduces to a gradient descent step after finitely many iterations, thus the proximal operator may be ignored. The proof follows \cite[Lemma 5.3]{Hale:2008} closely. 
  
  \begin{proposition}
  	\label{lemmings}\label{prop_lemmings}
  	For Prob.~(\ref{prob:sparse}), suppose Assumption 1 holds and $X^*$ is nonempty, thus there exists $h^*\in\mathbb{R}^n$ satisfying the conditions of Lemma \ref{thm:constGrad}. Let $E\triangleq\{i:|h_i^*|=\rho\}$.  Assume $\alpha_k=\beta_k$ for all $k$ and let $\{x^k,y^k\}_{k\in\mathbb{N}}$ be the iterates of (\ref{IFBS_yupdate})--(\ref{IFBS_xupdate}). Assume either:
  	\begin{enumerate}
  		\item $\{\lambda_k\}_{k\in\mathbb{N}}$ is nondecreasing, $0<\lambda_k\leq 1/L$ and $0\leq\alpha_k\leq\overline{\alpha}<1$ for all $k\in\mathbb{N}$, or
  		\item $\{\lambda_k,\alpha_k\}_{k\in\mathbb{N}}$ satisfy (\ref{chamchoice}).
  	\end{enumerate}
  	Let $x^* = \lim_{k\to\infty}x^k$ which exists by Theorem \ref{thm:fista_lyap} . Then for all but finitely many $k$, 
  	\begin{eqnarray}
  	x_i^{k} = y_i^{k}-\lambda_{k-1}(\nabla f(y^{k})_i - h_i^*),\quad\forall i\in \supp(x^*).\label{some}
  	\end{eqnarray}	
  \end{proposition}
  \emph{Proof.}
  See Sec.~\ref{proofOfTheorem55}.

  Under strict complementarity ($E=\supp(x^*)$) we will refer to the regime where (\ref{some}) is satisfied and $x_{E^c}^k=0$ as ``the large-$k$ regime" throughout the remainder of the paper. We refer to the regime where these conditions are not satisfied as ``the small-$k$ regime". 
  
  Now we consider a simple fixed parameter choice for Prob.~(\ref{prob:l1LS}). Under the strict complementarity condition, we can prove local linear convergence for any fixed choice of the inertia parameter in $[0,1)$ and the stepsize in $(0,1/L]$. The analysis turns out to be fairly elementary in this case since for this problem once in the large-$k$ regime the iterations form a simple 2\textsuperscript{nd} order linear homogeneous recursion which has been studied before, for example in \cite{o2012adaptive}. Note that we do not require the function $\phi$ defined in (\ref{eq:phidef}) to be strongly convex nor the minimizer to be unique. 
  
  \begin{theorem}
  	For Prob.~(\ref{prob:l1LS}) there exists $h^*\in\mathbb{R}^n$ satisfying the conditions of Lemma \ref{thm:constGrad}. Let $E\triangleq\{i:|h_i^*|=\rho\}$. Let $\beta_k=\alpha_k=\alpha\in[0,1)$ and $\lambda_k=1/L$ for all $k\in\mathbb{N}$. Let the iterates of (\ref{IFBS_yupdate})--(\ref{IFBS_xupdate}) be $\{x^k\}_{k\in\mathbb{N}}$ and $\lim_{k\to\infty} x^k = x^*$, which exists by Theorem \ref{thm:fista_lyap}. Suppose $E=\supp(x^*)$ (i.e. strict complementarity holds). Then $x^k$ achieves local Q-linear convergence. In particular there exists $K>0$, $C>0$, and $q\in(0,1)$ such that $\|x^k-x^*\|=C q^k$ for all $k>K$. Let $\hat{l}_E$ be the smallest eigenvalue of $A_E^\top A_E$. If $\hat{l}_E>0$, $0<\mu\leq\hat{l}_E$ and 
  	$
  	\alpha = (1-\sqrt{\mu/L})/(1+\sqrt{\mu/L})\label{param_choice_nonunique},
  	$ 
  	then $q=\left(1-\sqrt{\mu/L}\right)^{1/2}$. If $\hat{l}_E=0$ then $q\leq\alpha$.  Finally $F(x^k)$ converges to $F^*$ with rate $q^2$. 
  	
  	\label{cor:furtherconv} \label{thm:furtherconv} 
  \end{theorem} 
  \emph{Proof.}
  See Sec.~\ref{proofOfTheorem55}.

  For simplicity we prove the result only for $\lambda_k=1/L$ but the case $\lambda_k=\lambda\in(0,1/L]$ can also be shown. We stress that in practice the quantities $\hat{l}_E$ and $L$ are typically not known. 
  In the next section we show that a simple adaptive restart scheme can be incorporated into FISTA-CD to create a scheme which obtains the optimal iteration complexity without needing knowledge of $\hat{l}_E$.  
  
  
  
  \subsection{Asymptotic Behavior of FISTA-CD}\label{sec:asFISTA}
  We now ask, what is the convergence behavior of FISTA-CD in the large-$k$ regime? For Prob.~(\ref{prob:l1LS}) we see that once (\ref{some}) holds, the iterates are in the form of an inhomogeneous 2\textsuperscript{nd} order linear recurrence which has been studied previously in \cite{o2012adaptive} and \cite[\S 5--6]{tao2015local}. 
  It is difficult to analyze this recursion because $\alpha_k$ changes at each iteration and to do so rigorously requires a subtle argument following the one presented in \cite[\S 5--6]{tao2015local}, which is beyond the scope of this paper. 
  A simpler route to understanding the behavior is to use the homogeneous approximation of \cite[\S 4]{o2012adaptive} which sets $\alpha_k$ fixed and ``close" to $1$. This approximation implies that under strict complementarity, once in the large-$k$ regime and for $\alpha_k$ sufficiently close to $1$ (recall $\alpha_k\to 1$ for this parameter choice), FISTA-CD will exhibit nonmonotone oscillatory behavior in the objective function values with suboptimal Q-linear rate:
  \begin{eqnarray}
  \exists K,C>0: F(x^k)-F^*=C\left(\left(1-\lambda\hat{l}_E\right)^{k}\right),\quad\forall k>K\label{underdamp1},
  \end{eqnarray}
  where $\hat{l}_E$ is defined in Theorem \ref{thm:furtherconv}. This is the same as the convergence rate achieved by FBS (I-FBS with $\alpha_k=0$ and $\lambda_k=1/L$ for all $k\in\mathbb{N}$, although a slightly better rate can be achieved with $\lambda_k=2/(\hat{l}_E+L)$ which nevertheless has the same iteration complexity \cite{Hale:2008}). 

  For strongly convex quadratic minimization problems, \cite{o2012adaptive} suggested restarting the inertia sequence of Nesterov's method whenever a certain restart condition is observed. By applying the homogeneous approximation of \cite{o2012adaptive} to FISTA-CD once in the large-$k$ regime, we argue that we can improve the asymptotic convergence rate by incorporating such a restart technique. Thus even though the overall problem is nonsmooth and in general not strongly convex, restarting can improve the convergence  properties of FISTA-CD. 
  Restart schemes such as the ``speed restart" scheme \cite{su2014differential}, the ``gradient restart" scheme \cite{o2012adaptive}, the ``objective function" scheme \cite{o2012adaptive}, or the more conservative restart scheme of \cite{monteiro2012adaptive} could be incorporated into FISTA-CD. 
  For simplicity we elaborate only the objective function restart scheme of \cite{o2012adaptive} and we call the new method FISTA-CD-RE (``FISTA-CD with restart"). The idea is as follows. Whenever we observe $F(x^{k+1})>F(x^k)$, set the iteration counter $k$ in (\ref{chamchoice}) to $1$, and set $x^0=x^k$ and $x^1=x^k$. In other words restart FISTA-CD at the current point.  We refer the reader to \cite{o2012adaptive} for full details and analysis which can be applied to our situation in the large-$k$ regime (under strict complementarity). The homogeneous approximation of \cite{o2012adaptive} suggests FISTA-CD will have the optimal iteration complexity
  \begin{eqnarray*}
  \text{iter. comp.}= \Omega\lbr\sqrt{\frac{L}{\hat{l}_E}}\log1/\epsilon\rbr
  \end{eqnarray*}
  and rate
  \begin{eqnarray}
  F(x^k) - F^* = C'\left(1-\sqrt{\hat{l}_E/L}\right)^k.
  \label{eq_itercompFISTACD}
  \end{eqnarray}
  Remarkably it achieves this iteration complexity without knowledge of $\hat{l}_E$, the local strong convexity parameter.  Thus we do not need to know $\hat{l}_E$ in order to achieve the optimal convergence rate given in Theorem \ref{thm:furtherconv} with $\mu=\hat{l}_E$.  The method will also have $O(1/k^2)$ convergence rate while no restarts occur. It is also straightforward to incorporate a backtracking line search into FISTA-CD-RE, such as the one described in \cite[p.~194]{Beck:2009}, so that the method does not require $L$.
  
  We stress that the convergence rates given in (\ref{underdamp1}) and (\ref{eq_itercompFISTACD}) can be proved rigorously using arguments in the spirit of \cite[\S 5--6]{tao2015local}. However the analysis is beyond the scope of this paper. Our contribution is to show that, for all but finitely many iterations, FISTA-CD reduces to a form that has been studied previously in \cite{tao2015local,o2012adaptive}, from which convergence behavior can be extracted.

  \section{Numerical Simulations}
  \label{sec:sims}
  We now provide a small synthetic experiment to corroborate the theoretical findings of this paper. 
  \subsection{Experiment Details}
\begin{figure}
  	\centering
  	\includegraphics[width=4in]{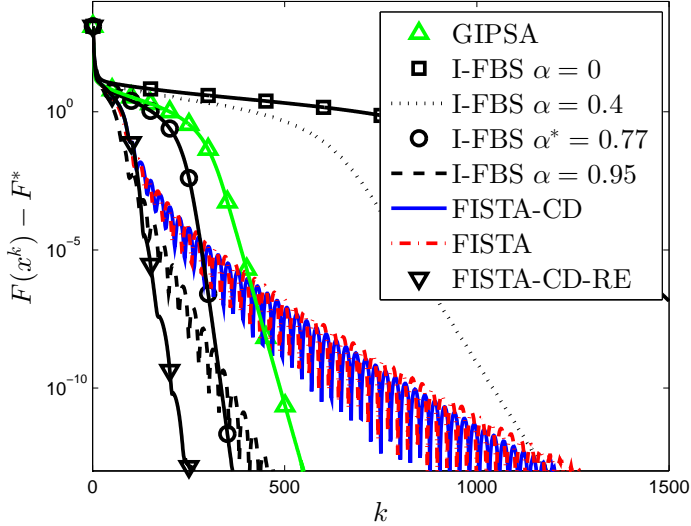}
  	\caption{Simulation results: showing $F(x^k)-F^*$ versus iteration $k$ for Experiment 1.}
  	\label{fig:sim1}
  \end{figure}

  We consider a randomly generated instance of Prob. (\ref{prob:l1LS}). The parameters of the experiment are $n=2000$, $m=1000$ and $\rho=0.1$. The entries of $A$ are drawn i.i.d. from the normal distribution with mean $0$ and variance $0.01$. The vector $b$ is given by $Ax_0$, where $x_0$ has $260$ nonzero entries generated i.i.d. from the $0$-mean unit variance normal distribution, and support set chosen uniformly at random. Recall that $l_E$ denotes the smallest nonzero eigenvalue of $A_E^\top A_E$ where $E$ is defined in Sec.~\ref{props_of_l1}. Note that for such a randomly generated problem where the entries of $A$ are drawn from a continuous probability distribution, $l_E>0$, and thus the solution is unique, with probability $1$ \cite{tibshirani2013lasso}. 
  We run (\ref{IFBS_yupdate})--(\ref{IFBS_xupdate}) with several choices for the parameters. For GIPSA: we choose $\beta_k=0.6$, $\alpha_k=0.42$ and $\lambda_k=1.39/L$ for all $k$. For I-FBS, $\lambda_k=1/L$ and $\beta_k=\alpha_k=\alpha\in\{0,0.4,\alpha^*,0.95\}$ where $\alpha^*$ is the locally optimal choice from Thm. \ref{cor:furtherconv}: $(1-\sqrt{l_E/L})/(1+\sqrt{l_E/L})$. The Lipschitz constant $L$ is the largest eigenvalue of $A^\top A$ and is computed via the SVD. We estimate $E$ via the interior point solver of \cite{kim2007interior} which we use to find an approximate solution $x^*$ such that the relative objective function error is no greater than $10^{-6}$. We then compute $h^*=\nabla f(x^*)$, and estimate $E$ as the set of all $i$ such that $\rho-|h^*_i|$ is smaller than $10^{-4}$. We then use the SVD of $A_E$ to estimate $l_E$. Using this approach, $\alpha^*$ is estimated as $0.77$ for this experiment. 
  Note that this is obviously not a practical method for estimating the optimal inertia parameter. The purpose of this experiment is simply to test the theoretical findings of Sections \ref{sec_mainResults}. In fact this experiment demonstrates that FISTA-CD-RE has the same asymptotic convergence rate as I-FBS with the optimal inertia parameter yet does not need to estimate $l_E$. We run FISTA, which is parameter choice \cite[Eq. (4.2)--(4.3)]{Beck:2009} with $\alpha_k=\beta_k$. We also run FISTA-CD which is (\ref{IFBS_yupdate})--(\ref{IFBS_xupdate}) with the parameter choice given in (\ref{chamchoice}) with $\lambda=1/L$ and $a=2.1$. We run FISTA-CD-RE with the same values for $\lambda$ and $a$ as FISTA-CD.  All algorithms are initialized to $x^1=x^{0}=0$. 
  The results are shown in Fig. \ref{fig:sim1} where we plot $F(x^k)-F^*$  versus $k$. Note the $y$-axis is logarithmic.
  
  \label{secblah}

  \subsection{Repeated Trials}
  We repeat this experiment 1000 times with different randomly drawn $A$ and $x_0$ from the distributions described above. For each trial we record the number of iterations after which the the relative error remains below tol, i.e. $k:(F(x^j)-F^*)/F^* \leq \text{tol},\forall j\geq k$.\footnote{$F^*$ is approximated by the smallest objective function value among all tested algorithms after $1500$ iterations} The average of this number across the $1000$ trials is given in Table \ref{tableux} for $\text{tol}\in\{10^{-2},10^{-6}\}$ and all algorithms. 
  \label{sec_Repeated_trials}
  
  \begin{table}
  	\caption{Results for repeated trials (Sec.~\ref{sec_Repeated_trials})}\label{tableux}
  	\begin{tabular}{ ||p{3cm}||p{3.5cm}|p{3.5cm}||  }
  		\hline
  		Algorithm& Average \# iterations to rel. err. $10^{-2}$ (1000 trials)& Average \# iterations to rel. err. $10^{-6}$ (1000 trials)\\
  		\hline
  		GIPSA & 260 & 368 \\
  		\hline
  		I-FBS ($\alpha=0$)   & 901    & 1287\\
  		I-FBS ($\alpha=0.4$)&   540  & 775   \\
  		I-FBS ($\alpha=\alpha^*$)& 210 & 286\\
  		I-FBS ($\alpha=0.95$)& 68 & 171 \\
  		\hline
  		FISTA & 84 & 282 \\
  		FISTA-CD&   85  & 280\\
  		\hline
  		FISTA-CD-RE& 85  & 137  \\
  		\hline
  	\end{tabular}
  	
  \end{table}
  
  \subsection{Observations}
  First let's look at Fig. \ref{fig:sim1}. Although the figure shows objective function values, since the minimizer is unique, convergence of the iterates is implied, which corroborates Theorem \ref{thm:fista_lyap} (not including FISTA, for which convergence of the iterates is an open problem). All tested parameter choices for I-FBS transition from a manifold identification period to a local linear convergence period, corroborating Theorem \ref{thm:finite}. Furthermore adding inertia does improve the asymptotic rate and using $\alpha^*$ achieves the best asymptotic rate. However FISTA-CD-RE essentially achieves the same asymptotic rate despite not knowing or estimating $l_E$. The upper bound for the asymptotic convergence rate of I-FBS with inertial parameter $\alpha^*$ is computed using (\ref{bnn}) to be $0.89$, which compares with an empirically determined rate of $0.83$. However the fixed choice $\alpha_k=\alpha^*$ is outperformed by the larger choice $\alpha=0.95$, and FISTA, FISTA-CD and FISTA-CD-RE in the small-$k$ regime (i.e. before linear convergence commences). In the large-$k$ regime FISTA-CD exhibits nonmonotone oscillatory behavior and suboptimal asymptotic convergence as predicted in Sec.~\ref{sec:asFISTA}. 
  
  Now we look at Table \ref{tableux}. For a low accuracy solution (i.e. rel. err. less than $10^{-2}$), I-FBS with $\alpha=0.95$ performs best and there is no difference between FISTA, FISTA-CD and FISTA-CD-RE. FISTA-CD and FISTA-CD-RE are identical because a restart had not yet occurred in any of the 1000 trials. However for a high accuracy solution (i.e. rel. err. less than $10^{-6}$), FISTA-CD-RE outperforms all other algorithms. It requires on average fewer than half as many iterations as FISTA or FISTA-CD at essentially the same per-iteration cost \footnote{Despite having an additional function evaluation per iteration, FISTA-CD-RE only requires one matrix multiply per iteration, which is the same as FISTA-CD and FISTA since the matrix multiply is the dominant cost}. Note that despite the optimal asymptotic convergence rate of the fixed choice $\alpha^*$, this choice is outperformed when looking for a high accuracy solution. This is because its small-$k$ regime performance is poor. The strong performance of I-FBS with $\alpha=0.95$ in the small-$k$ regime is interesting and we cannot explain it with the existing theory. However for such large values of the inertia parameter we expect the performance to be approximately similar to FISTA and its variants.

  \section{Proof of Lemma \ref{thm:constGrad}}\label{proofOfLemConstGrad1}
  \label{proofOfLemConstGrad}
  We commence by proving that the gradient with respect to $f$ is constant at all optimal points. 
  The proof follows by considering \cite[Corollary 26.3(vii)]{bauschke2011convex}. Note that condition (a) of this Corollary holds trivially because $\text{dom} f = \mathcal{H}$ and $\text{dom}\,\partial g\subseteq\text{dom}\,g$ is nonempty. Now statement (vii) of Corollary 26.3 states the following. Given $x\in X^*$ 
  \begin{eqnarray}
  \langle x-y,\nabla f(x)\rangle + g(x)\leq g(y)\quad\forall y\in\mathcal{H}.\label{lates}
  \end{eqnarray} 
  Consider $x_1,x_2\in X^*$, then (\ref{lates}) implies $\langle x_1-x_2,\nabla f(x_1)\rangle + g(x_1)\leq g(x_2)$ and  $\langle x_2-x_1,\nabla f(x_2)\rangle + g(x_2)\leq g(x_1)$. Adding these two together yields
  $$
  \langle\nabla f(x_1)-\nabla f(x_2),x_1-x_2\rangle\leq 0.
  $$
  From this point on the proof is identical to \cite[Prop.~26.10]{bauschke2011convex}, which implies $\nabla f(x_1)=\nabla f(x_2)\triangleq h^*$.
  The rest of the Lemma follows by examining the structure of the optimality condition (\ref{optcond}) for the special case of Prob.~(\ref{prob:sparse}). We refer the reader to \cite[Thm 2.1 (ii) and (iii)]{Hale:2008}.\qed
  

  \section{Proof of Theorem \ref{thm:fista_lyap}}\label{proofOfThm41}
  Before proving the theorem, we give three lemmas, beginning with the celebrated lemma due to Opial.
  \begin{lemma}[\cite{opial1967weak}, Opial's lemma]
  	\label{opiallemma}
  	Suppose $\{x^k\}$ is a sequence in $\mathcal{H}$ and $S\subset\mathcal{H}$ is a nonempty set such that:
  	\begin{enumerate}
  		\item $\lim_{k\to\infty}\|x^k-x^*\|$ exists for every $x^*\in S$,
  		\item Every weakly convergent subsequence of $\{x^k\}_{k\in\mathbb{N}}$ weakly converges to some $x^*\in S$.
  	\end{enumerate}
  	Then there exists $\hat{x}\in S$ such that $x^k\rightharpoonup\hat{x}$. 
  \end{lemma}
  
  It is trivial to verify that the second condition of Opial's lemma holds for GIPSA, so long as $x^k - x^{k-1}\to 0$. We do this in the following Lemma.
  \begin{lemma} 
  	\label{fixedpointlemma}
  	For Prob.~(\ref{prob:1}) suppose Assumption 1 holds and $X^*$ is nonempty. Let $\{x^k\}_{k\in\mathbb{N}}$ be the sequence generated by (\ref{IFBS_yupdate})--(\ref{IFBS_xupdate}). Suppose $x^{v_k}\rightharpoonup x$ for some subsequence $\{v_k\}_{k\in\mathbb{N}}\subseteq \mathbb{N}$, and $x^k-x^{k-1} \to 0$.  Then $x\in X^*$. 
  \end{lemma}
  \emph{Proof.}
  	The proof follows the techniques of \cite[Theorem 2.1]{MoudafiOliny}. Thanks to (\ref{IFBS_yupdate}) and the assumption that $x^k-x^{k-1}\to 0$, we know that $y^{k+1}-x^k\to 0$ and thus $x^k-y^k \to 0$. Similarly by (\ref{IFBS_zupdate}) we see that $z^{k+1}-x^k\to 0$ and therefore $z^k - x^k \to 0 $. Now by (\ref{implicitupdate}) 	
  	\begin{eqnarray}\label{hi}
  	-\frac{1}{\lambda_{v_k - 1}}(x^{v_k}-y^{v_k}) + \nabla f(x^{v_k})-\nabla f(z^{v_k})&\in & \{\nabla f(x^{v_k})\}+\partial g(x^{v_k}).
  	\end{eqnarray}
  	Now passing to the limit $v_k\to\infty$, using the fact that $\nabla f$ is Lipschitz continuous, and \cite[Proposition 3.4(b)]{burachik1999varepsilon}, we infer that $0\in\partial g(x)+\{\nabla f(x)\}$, therefore $x\in X^*$ by optimality condition (\ref{optcond}).
  	\qed 
  
  The final Lemma is standard in the analysis of inertial methods. 
  \begin{lemma}
  	Let $\{\varphi_k,\delta_k,\sigma_k\}_{k\in\mathbb{N}}\subset \mathbb{R}_+$. If $\varphi_{k+1}-\varphi_k\leq \sigma_k(\varphi_k - \varphi_{k-1}) + \delta_k$ for all $k$ where $\sigma_k\leq\overline{\sigma}<1$ and $\sum_{k\in\mathbb{N}}\delta_k <\infty$, then $\lim_{k\to\infty}\varphi_k$ exists. \label{theta_lem}
  \end{lemma}
  
  \emph{Proof.} We refer to \cite[Thm 3.1]{alvarez2000minimizing}.\qed
  
  We now turn our attention to Theorem \ref{thm:fista_lyap}. We prove statement (\ref{state:finitesum}) by using the multistep Lyapunov function from  \cite{alvarez2000minimizing} which is shown to be nonincreasing. The proof of (\ref{criticality}) is trivial. Finally to prove (\ref{state:weak}) we use Lemma \ref{theta_lem} to prove the first condition of Opial's lemma holds (the second condition of Opial's lemma holds by Lemma \ref{fixedpointlemma}).
  
  \subsubsection*{Proof of Theorem \ref{thm:fista_lyap} statement (\ref{state:finitesum})}
  Recall the notation: $\Delta_k\triangleq x^k - x^{k-1}$. Define the energy function:
  $
  E_{k} \triangleq F(x^k) + \frac{\beta_k}{2\lambda_k}\|\Delta_k\|^2.
  $
  We will show that $E_k$ is nonincreasing. Using the fact that $\lambda_k$ is nondecreasing, then (\ref{eq:ineq2}) and (\ref{eq:ineq1}),  we write
  \begin{eqnarray}
  E_{k+1}-E_k
  &\leq&
  F(x^{k+1})-F(x^k)
  +\frac{1}{2\lambda_k}\left(\beta_{k+1}\|\Delta_{k+1}\|^2-\beta_{k}\|\Delta_k\|^2\right)
  \nonumber\\
  &\leq&
  \left\langle\nabla f(z^{k+1})+v,\Delta_{k+1}\right\rangle
  +
  \frac{L}{2}\|z^{k+1} - x^{k+1}\|^2
  \nonumber\\&&\label{mays}
  +
  \frac{1}{2\lambda_k}\left(\beta_{k+1}\|\Delta_{k+1}\|^2-\beta_{k}\|\Delta_k\|^2\right),\quad\forall v\in\partial g(x^{k+1}).
  \end{eqnarray}
  Note that by the definition of the $\prox$ operator, $x^{k+1}\in\text{dom}\, \partial g$ which is nonempty by Assumption 1. Using (\ref{implicitupdate}) in (\ref{mays}) implies
  \begin{eqnarray}
  E_{k+1}-E_k&\leq&
  \frac{1}{\lambda_k}\left\langle y^{k+1} - x^{k+1},\Delta_{k+1}\right\rangle
  +\frac{L}{2}\|z^{k+1} - x^{k+1}\|^2
  \nonumber\\&&\label{latest}
  +
  \frac{1}{2\lambda_k}\left(\beta_{k+1}\|\Delta_{k+1}\|^2-\beta_{k}\|\Delta_k\|^2\right).
  \end{eqnarray}
  Now using (\ref{IFBS_yupdate}) and (\ref{IFBS_zupdate}) we derive:
  \begin{eqnarray}
  y^{k+1} - x^{k+1} = \beta_k\Delta_k - \Delta_{k+1}
  \,\,\text{ and }\,\,
  z^{k+1} - x^{k+1} = \alpha_k\Delta_k - \Delta_{k+1}\label{some_facts}.
  \end{eqnarray}
  Substituting (\ref{some_facts}) into (\ref{latest}) yields
  \begin{eqnarray}
  E_{k+1}-E_k\nonumber
  &\leq&
  \left(\frac{\beta_k-\alpha_k\lambda_k L}{\lambda_k}\right)\langle\Delta_{k+1},\Delta_k\rangle
  +\frac{\alpha_k^2 \lambda_k L - \beta_k}{2\lambda_k}\|\Delta_k\|^2
  \nonumber\\&&
  +\left(\frac{\beta_{k+1}}{2\lambda_k} 
  + \frac{L}{2}-\frac{1}{\lambda_k}\right)\|\Delta_{k+1}\|^2
  \nonumber\\
  &=&
  -\left(\frac{\beta_k-\alpha_k\lambda_k L}{2\lambda_k}\right)\|\Delta_{k+1}-\Delta_k\|^2
  -\frac{\alpha_k(1-\alpha_k) L}{2}\|\Delta_k\|^2
 \nonumber \\\label{aga}
  &&
  -
  \frac{2-\lambda_k L(1-\alpha_k) - \beta_k - \beta_{k+1}}{2\lambda_k}
  \|\Delta_{k+1}\|^2.
  \end{eqnarray}
  Now (\ref{param_conditions}) implies that the coefficients of $\|\Delta_{k+1}-\Delta_k\|^2$, $\|\Delta_k\|^2$ and $\|\Delta_{k+1}\|^2$ are nonpositive.
  Furthermore, from condition (\ref{param_conditions}) we see that
  \begin{eqnarray*}
  	\frac{2-\lambda_k L(1-\alpha_k)-\beta_k - \beta_{k+1}}{2\lambda_k}\geq \frac{\epsilon}{2\lambda_k}> \frac{\epsilon L}{4}>0.
  \end{eqnarray*}
  Therefore telescoping (\ref{aga}) implies 
  \begin{eqnarray*}
  	\frac{\epsilon L}{4}\sum_{k=1}^M \|\Delta_{k+1}\|^2 < E_1 - E_{M+1} <\infty,\quad\forall M\in\mathbb{N}.
  \end{eqnarray*}
  Thus statement (i) is proven. 
  
  
  \subsubsection*{Proof of Theorem \ref{thm:fista_lyap} statement (\ref{criticality})}
  
  Consider (\ref{hi}) with the subsequence chosen as $v_k=k$. Clearly this implies statement (\ref{criticality}) since the left-hand side of (\ref{hi}) goes to $0$ as $k$ goes to $\infty$. 
  \subsubsection*{Proof of Theorem \ref{thm:fista_lyap} statement (\ref{state:weak})}
  Assume $X^*$ is nonempty and $x^{v_k}$ is a subsequence which weakly converges to $x'$. We note that statement (i) implies $\Delta_k\to 0$, thus from Lemma \ref{fixedpointlemma}, $x'\in X^*$. Therefore the second condition of Opial's lemma is satisfied. 
  
  We now proceed to show the first condition of Opial's lemma, i.e. for any $x^*\in X^*$, the limit of $\{\|x^k-x^*\|\}_{k\in\mathbb{N}}$ exists. The key will be to derive a recursion in the form of Lemma \ref{theta_lem}. This part of the proof has been adapted from \cite{MoudafiOliny} which studies the special case where $\alpha_k=0$ for all $k$. 
  Fix $x^*\in X^*$ (which is nonempty by assumption) and let $\varphi_k\triangleq \frac{1}{2}\|x^k - x^*\|^2$. Now using (\ref{biglem}) we see that
  $$
  \langle x^{k+1}-x^k, x^*-x^{k+1}\rangle = \varphi_k - \varphi_{k+1} - \frac{1}{2}\|\Delta_{k+1}\|^2.\label{biglem_spec}
  $$
  Combining this with (\ref{IFBS_yupdate}) yields
  \begin{eqnarray}
  \varphi_k - \varphi_{k+1}
  &=&
  \frac{1}{2}\|\Delta_{k+1}\|^2 
  + 
  \langle x^{k+1}-y^{k+1},x^*-x^{k+1}\rangle
  \nonumber\\&&\label{mmo}
  +\,
  \beta_k\langle x^k - x^{k-1},x^*-x^{k+1}\rangle.
  \end{eqnarray}
  Now by (\ref{implicitupdate}):
  $$
  -\left(x^{k+1}-y^{k+1} +\lambda_k \nabla f(z^{k+1})\right)\in\lambda_k \partial g(x^{k+1}).
  $$
  On the other hand by optimality condition (\ref{optcond}):
  $
  -\lambda_k \nabla f(x^*)\in\lambda_k \partial g(x^*).
  $
  Using these facts and (\ref{montone_op}) gives
  \begin{eqnarray} 
  \langle x^{k+1}-y^{k+1} +\lambda_k(\nabla f(z^{k+1})-\nabla f(x^*)),x^*-x^{k+1}\rangle \geq 0
  \nonumber
  \end{eqnarray}
  which implies
  \begin{eqnarray} 
  \label{but}
  \langle x^{k+1}-y^{k+1},x^*-x^{k+1}\rangle \geq
  \lambda_k\langle \nabla f(z^{k+1})-\nabla f(x^*) ,x^{k+1}-x^*\rangle.
  \end{eqnarray}
  Substituting (\ref{but}) into (\ref{mmo}) yields 
  \begin{eqnarray}
  \varphi_{k+1} - \varphi_{k}
  &\leq&
  -\frac{1}{2}\|\Delta_{k+1}\|^2 
  -
  \lambda_k\langle \nabla f(z^{k+1})-\nabla f(x^*) ,x^{k+1}-x^*\rangle
  \nonumber\\&&\label{xax}
  +\,
  \beta_k\langle x^k - x^{k-1},x^{k+1}-x^*\rangle.
  \end{eqnarray}
  Now using (\ref{biglem}) again
  \begin{eqnarray}
 \langle x^k - x^{k-1},x^{k+1}-x^*\rangle
 =&&
  \varphi_k - \varphi_{k-1} +\frac{1}{2}\|\Delta_k\|^2 
  \nonumber\\\label{wow}
  &&+\,\langle x^k - x^{k-1},x^{k+1}-x^k\rangle.
  \end{eqnarray}
  On the other hand using (\ref{def_coercive})
  \begin{eqnarray}
  \langle \nabla f(z^{k+1})-\nabla f(x^*),x^{k+1}-x^*\rangle
  &\geq&
  \frac{1}{L}
  \|\nabla f(z^{k+1})-\nabla f(x^*)\|^2 
  \nonumber\\
  &&
  + \left\langle \nabla f(z^{k+1}-\nabla f(x^*),x^{k+1}-z^{k+1}\right\rangle
  \nonumber\\
  &\geq&
  -\frac{L}{4}\|x^{k+1}-z^{k+1}\|^2.\label{hubba}
  \end{eqnarray}
  Therefore by substituting (\ref{wow}) and (\ref{hubba}) into (\ref{xax}) and using (\ref{some_facts}), we get
  \begin{eqnarray}
  \varphi_{k+1}-\varphi_k-\beta_k(\varphi_k - \varphi_{k-1})
  \leq&&
  -\frac{\zeta_k}{2}\|\Delta_{k+1}\|^2
  +c^k_1\|\Delta_k\|^2
  \nonumber\\\label{late}
  &&
  +
  c^k_2\langle\Delta_k,\Delta_{k+1}\rangle,
  \end{eqnarray}
  where
  $
  \zeta_k = (1-\lambda_k L/2),
  $
  $
  c_1^k = (\beta_k/2 +\alpha_k^2\lambda_k L/4),
  $
  and
  $
  c^k_2 = \beta_k - \alpha_k\lambda_k L/2. 
  $
  Note that (\ref{param_conditions}) implies $\zeta_k\geq\gamma/2>0$, $c_1^k\in[0,1)$ and $|c_2^k|<1$. Now if we let $\theta_k \triangleq \varphi_k - \varphi_{k-1}$ then (\ref{late}) implies
  \begin{eqnarray*}
  	\theta_{k+1}-\beta_k\theta_k
  	&\leq&
  	-\frac{\zeta_k}{2}
  	\left\|
  	\Delta_{k+1} - \frac{c^k_2}{\zeta_k}\Delta_k
  	\right\|^2
  	+
  	\left(c_1^k+\frac{(c_2^k)^2}{2\zeta_k}\right)\|\Delta_{k}\|^2\leq \delta_k,
  \end{eqnarray*}
  with $\delta_k \triangleq \left(1+1/\gamma\right)\|\Delta_{k}\|^2$.
  Note that (\ref{state:finitesum}) of this Theorem implies $\sum_{k\in\mathbb{N}}\delta_k <\infty$. Now since $\beta_k\leq\overline{\beta}<1$, we can apply Lemma \ref{theta_lem}, which implies $\lim_{k\to\infty}\|x^k-x^*\|$ exists for any $x^*\in X^*$. Therefore both conditions of Opial's lemma hold and $\{x^k\}_{k\in\mathbb{N}}$ converges weakly to some minimizer $\hat{x}$.
  Now repeating (\ref{IFBS_yupdate}): $y^{k+1}=x^k+\beta_k(x^k-x^{k-1})$ for all $k\in \mathbb{N}$. Therefore for any $h\in\mathcal{H}$, 
  $\langle h,y^{k+1}\rangle=\langle h,x^k\rangle +\langle h,\beta_k(x^k-x^{k-1})\rangle\to\langle h,\hat{x}\rangle
  $, which proves $y^k\rightharpoonup\hat{x}$. In exactly the same way we can show $z^k\rightharpoonup \hat{x}$ using (\ref{IFBS_zupdate}).\qed
  
  \section{Proof of Theorem \ref{thm:finite}}
  \label{proofOfTheroem5}
  Before proving the theorem, we need the following three lemmas. The first details the contractive properties of the soft-thresholding operator.  
  
  \begin{lemma}[\cite{Hale:2008}, Lemma 3.2]
  	Fix any $a$ and $b$ in $\mathbb{R}$, and $\nu\geq 0$:
  	\begin{enumerate}[(i)]
  		\item \cite[Lemma 3.2 (3.7)]{Hale:2008} The function $S_\nu$ defined in (\ref{def_softThresh}) is nonexpansive. That is,\label{eq:lemmaSv1}
  		$
  		|S_\nu(a)-S_\nu(b)|\leq|a-b|.
  		$
  		\item \cite[Lemma 3.2 statement (5)]{Hale:2008} If $|b|\geq\nu$ and $\sgn(a) \neq \sgn(b)$ then \label{eq:lemmaSv2}
  		$
  		|S_{\nu}(a)-S_{\nu}(b)|\leq|a-b|-\nu.
  		$
  		\item \cite[Lemma 3.2 statement (6)]{Hale:2008} If $S_{\nu}(a)\neq 0=S_{\nu}(b)$ then 
  		$
  		|S_{\nu}(a)-S_\nu(b)|\leq|a-b|-(\nu-|b|).
  		$\label{eq:lemmaSv3}
  	\end{enumerate}
  	\label{lemma:Sv}
  \end{lemma}
  
  Next we derive some technical properties of the solution set for Prob.~(\ref{prob:sparse}).
  \begin{lemma}\label{lemmaeq12}\label{callback2}
  	For Prob.~(\ref{prob:sparse}) suppose Assumption 1 holds and $X^*$ is nonempty, $x^*\in X^*$ and $\lambda>0$. Then there exists a vector $h^*\in\mathbb{R}^n$ satisfying the conditions of Lemma \ref{thm:constGrad}. Furthermore
  	\begin{eqnarray}
  	|x_i^*-\lambda h_i^*|\geq\rho\lambda,
  	\label{lemmaeq1}
  	\text{ and }\,\,
  	\sgn(x_i^*-\lambda h_i^*) = -h_i^*/\rho,\quad \forall i:|h_i^*|=\rho.\label{callback}
  	\end{eqnarray}
  \end{lemma}
  \emph{Proof.}
  	Recall that $E\triangleq\{i:|h_i^*|=\rho\}$.	For $i\in\supp (x^*)$, (\ref{def_softThresh}) and (\ref{optcond_fixed}) imply
  	\begin{eqnarray}
  	0\neq x_i^* = \sgn\left(x_i^*-\lambda h^*_i\right)\left[|x_i^*-\lambda h^*_i|-\rho\lambda\right]_+.
  	\label{eq:fix_point}
  	\end{eqnarray}
  	Therefore $|x_i^*-\lambda h^*_i|>\rho\lambda$ for all $i\in\supp(x^*)$. On the other hand, if $i\in E\setminus \supp(x^*)$, then 
  	$
  	|x_i^*-\lambda h^*_i|=\lambda|h_i^*|=\rho\lambda.
  	$
  	Recall that $\supp(x^*)\subseteq E$. Therefore the first part of (\ref{lemmaeq1}) is proven. 
  	
  	Looking at (\ref{eq:fix_point}) it can be seen that
  	\begin{eqnarray}
  	\sgn(x^*_i)=\sgn(x_i^*-\lambda h^*_i),\quad\forall i\in \supp(x^*).
  	\label{eq:l1ls2}
  	\end{eqnarray}
  	Note by Lemma \ref{thm:constGrad}, if $i\in \supp(x^*)$, then $\sgn(x_i^*)=-h_i^*/\rho$. Else if $i\in E\setminus \supp(x^*)$ then
  	\begin{eqnarray}
  	\sgn(x_i^*-\lambda h^*_i)= \sgn(-\lambda h_i^*)
  	=-\sgn(h_i^*)
  	= -\frac{h_i^*}{\rho}.\label{eq:l1ls1}
  	\end{eqnarray}   
  	since $|h_i^*|=\rho$. Combining (\ref{eq:l1ls2}) and (\ref{eq:l1ls1}) yields the second part of (\ref{callback}).
  \qed	
  
  The final lemma before we proceed with the proof of Theorem \ref{thm:finite} is a crucial finite summability result.  
  \begin{lemma}
  	\label{lem_recurse}
  	For Prob.~(\ref{prob:sparse}) suppose Assumption 1 holds. Assume either
  	\begin{enumerate}
  		\item $\{\lambda_k\}_{k\in\mathbb{N}}$ is nondecreasing, $0<\lambda_k\leq1/L$, $\alpha_k=\beta_k$ and $0\leq\alpha_k\leq \overline{\alpha}<1$  for all $k\in\mathbb{N}$, or
  		\item $X^*$ is nonempty and $\{\alpha_k,\lambda_k\}_{k\in\mathbb{N}}$ are chosen according to FISTA-CD in (\ref{chamchoice}).
  	\end{enumerate}
  	Furthermore assume the iterates $\{x^k,y^k\}_{k\in\mathbb{N}}$ of (\ref{IFBS_yupdate})--(\ref{IFBS_xupdate}) satisfy, for all $k\in\mathbb{N}$:
  	\begin{eqnarray}
  	\quad\quad\|x^{k}-x\|^2 \leq \|y^{k} - x\|^2 - N_{k}\label{eq:lem_recurs}
  	\end{eqnarray}
  	for some $x\in \mathbb{R}^n$ and $\{N_k\}_{k\in\mathbb{N}}\subset\mathbb{R}_+$. Then:
  	$
  	\sum_{k=1}^\infty N_k <\infty.
  	$
  \end{lemma}
  
  \vspace{0.2cm}
  
  \emph{Proof.}
  	Substituting (\ref{IFBS_yupdate}) into (\ref{eq:lem_recurs}) yields
  	\begin{eqnarray}
  	\|x^{k+1}-x\|^2
  	&\leq&
  	\|x^k-x + \alpha_k\Delta_k\|^2 - N_{k+1}
  	\nonumber\\
  	&=&
  	\|x^k - x\|^2 + \alpha_k^2\|\Delta_k\|^2 + 2\alpha_k\langle x^k - x,\Delta_k\rangle - N_{k+1}\label{simple}.
  	\end{eqnarray}
  	Let $\varphi_k=\frac{1}{2}\|x^k - x\|^2$ and $\theta_k=\varphi_k - \varphi_{k-1}$.  Using (\ref{biglem}) we write
  	$
  	\langle x^k - x,\Delta_k\rangle = \varphi_k - \varphi_{k-1}+\frac{1}{2}\|\Delta_k\|^2.
  	$
  	Using this in (\ref{simple}) yields
  	\begin{eqnarray}
  	\theta_{k+1}\leq \alpha_k\theta_k + \delta_k - \frac{1}{2}N_{k+1}\label{asd1}
  	\end{eqnarray}
  	where $\delta_k=\frac{1}{2}\alpha_k(1+\alpha_k)\|\Delta_k\|^2$. Note that $0\leq\delta_k\leq\alpha_k\|\Delta_k\|^2\leq\|\Delta_k\|^2$. 
  	
  	We first prove the lemma for parameter choice 1. For this parameter choice by Theorem \ref{thm:fista_lyap}(\ref{state:finitesum}), $\sum_{k\in\mathbb{N}}\delta_k <\infty$. Let $\underline{\alpha}=\inf_k\alpha_k$ and note that $\underline{\alpha}\in[0,\overline{\alpha}]$. Thus using (\ref{asd1}):
  	\begin{eqnarray*}
  		\theta_{k+1}
  		&\leq&
  		\overline{\alpha}^k|\theta_1|+ \sum_{j=1}^k\overline{\alpha}^{k-j}\delta_j - \frac{1}{2}\sum_{j=1}^k\underline{\alpha}^{k-j}N_{j+1}.
  	\end{eqnarray*}
  	Therefore, for all $M\in\mathbb{N}$,
  	\begin{eqnarray}
  	\varphi_M
  	&=&
  	\varphi_0 + 
  	\sum_{k=1}^M\theta_{k}
  	\leq
  	\varphi_0 
  	+
  	\frac{1}{1-\overline{\alpha}}\left(|\theta_1|+\sum_{k=1}^{M-1}\delta_k\right)
  	-\frac{1}{2}\sum_{k=1}^{M-1} N_{k+1}
  	\label{eq:renew1}\nonumber
  	\\\nonumber
  	&\implies&
  	\sum_{k=1}^{M-1}N_{k+1}
  	\leq
  	2\varphi_0
  	+\frac{2}{1-\overline{\alpha}}
  	\left(|\theta_1|+\sum_{k=1}^{\infty}\delta_k\right)<\infty,\quad\forall M\in\mathbb{N}.
  	\end{eqnarray} 
  	
  	Now for parameter choice 2, we proceed as follows. Note that since $x^1=x^0$, $\theta_1=0$ for this parameter choice. From (\ref{asd1}), and $\delta_k\leq\alpha_k\|\Delta_k\|^2$, we infer (using the convention: $\prod_{j=a}^b\alpha_j = 1$ if $a>b$):
  	\begin{eqnarray*}
  		\theta_{k+1}
  		&\leq&
  		\left(\prod_{i=1}^k\alpha_i\right)\theta_1+
  		\sum_{j=1}^k\left(\prod_{l=j}^k\alpha_l\right)\|\Delta_j\|^2 
  		-\frac{1}{2}\sum_{j=1}^k\left(\prod_{l=j+1}^k\alpha_l\right)N_{j+1}
  		\\
  		&\leq &
  		\sum_{j=1}^k\left(\prod_{l=j}^k\alpha_l\right)\|\Delta_j\|^2
  		-\frac{1}{2}\sum_{j=1}^k \alpha_2^{k-j}N_{j+1}
  	\end{eqnarray*}
  	where we have used the fact that $\alpha_2<\alpha_k$ for all $k>2$. Thus for all $M\in\mathbb{N}$
  	\begin{eqnarray}
  	\nonumber
  	\varphi_M
  	=
  	\varphi_0 + 
  	\sum_{k=1}^M\theta_{k}
  	&\leq&
  	\varphi_0 
  	+\sum_{k=1}^{M-1}
  	\sum_{j=1}^k\left(\prod_{l=j}^k\alpha_l\right)\|\Delta_j\|^2
  	-\frac{1}{2}\sum_{k=1}^{M-1} N_{k+1}\left(\sum_{j=0}^{M-1-k}\alpha_2^j\right)
  	\\
  	&\leq &\label{whaa1}\nonumber
  	\varphi_0  
  	+\sum_{k=1}^{\infty}
  	\sum_{j=1}^k\left(\prod_{l=j}^k\alpha_l\right)\|\Delta_j\|^2
  	-\frac{1}{2}\sum_{k=1}^{M-1} N_{k+1}.
  	\end{eqnarray}
  	Now by applying (\ref{abiggy}) of Lemma \ref{lemCham} and noting that $\varphi_M\geq 0$, we infer $\sum_{k=1}^\infty N_{k}<\infty$.  
  \qed
  
  We are now ready to prove Theorem \ref{thm:finite}. Note that parameter choice 1 satisfies the requirements of Corollary \ref{cor_IFBS_weak}. Furthermore, by assumption, $X^*$ is nonempty, thus all conclusions of Corollary \ref{cor_IFBS_weak} hold. For parameter choice 2 (FISTA-CD) we note that both conclusions of Lemma \ref{lemCham} hold. 
  
  Throughout the proof, fix an arbitrary $x^*\in X^*$. We will use the contractive properties of $S_\nu$ given in Lemma \ref{lemma:Sv} to construct a recursion in the form of (\ref{eq:lem_recurs}) of Lemma \ref{lem_recurse}. That lemma allows us to argue that the number of iterations such that (\ref{thm:finite:eq2})--(\ref{thm:finite:eq1}) do not hold is finite. Note that since $\alpha_k=\beta_k$ for the parameter choices in question (i.e. I-FBS and FISTA-CD), (\ref{IFBS_xupdate}) simplifies to
  \begin{eqnarray}\label{xupdate_special_ifbs}
  x^{k+1}=\prox_{\lambda_k g}\left(y^{k+1} - \lambda_k\nabla f(y^{k+1})\right)=T_{\lambda_k}(y^{k+1})
  \end{eqnarray}
  where $T_\lambda$ is the forward-backward operator discussed in Sec.~\ref{sec_assump_optimal}. 
  
  \subsubsection*{Proof of (\ref{thm:finite:eq2}) of Theorem \ref{thm:finite}} 
  Recall from Lemma \ref{thm:constGrad} there exists a vector $h^*$ such that $\nabla f(x^*) = h^*$ for all $x^*\in X^*$, and $\supp (x^*)\subseteq E$, where $E\triangleq\{i:|h_i^*|=\rho\}$. Fix $k\in\mathbb{N}$. 
  Now (\ref{xupdate_special_ifbs}) and optimality condition (\ref{optcond_fixed}) imply 
  \begin{eqnarray}
  |x_i^{k+1}-x_i^*|^2&=&\left|S_{\rho\lambda_k} (y_i^{k+1}-\lambda_k\nabla f(y^{k+1})_i)-S_{\rho\lambda_k} (x_i^*-\lambda_k h^*_i)\right|^2\label{clearing}
  \end{eqnarray}
  for all $i\in[n]$, using the notation $[n]\triangleq\{1,2,\ldots,n\}$.
  Consider the following condition:
  \begin{eqnarray}
  \hspace{0.5cm}\sgn\left(y^{k+1}_i-\lambda_k\nabla f(y^{k+1})_i\right)\neq \sgn(x_i^*-\lambda_k h^*_i)\quad\hbox{for some $i\in E$}.
  \label{eq:sgn}
  \end{eqnarray}
  (Note that $\sgn(x_i^*-\lambda_k h^*_i)=-h_i^*/\rho$ from Lemma \ref{callback2}). Now recall Lemma \ref{callback2} implies $|x_i^*-\lambda_k h^*_i|\geq \lambda_k\rho$ for all $i\in E$. Therefore we can apply Lemma \ref{lemma:Sv} (\ref{eq:lemmaSv2}) to (\ref{clearing}) to say the following. If (\ref{eq:sgn}) holds, then 
  \begin{eqnarray}
  \nonumber|x_i^{k+1}-x_i^*|^2
  &\leq&\left(|y_i^{k+1}-\lambda_k\nabla f(y^{k+1})_i-(x_i^*-\lambda_k h^*_i)|-\rho\lambda_k\right)^2\\
  &\leq&\left|y_i^{k+1}-\lambda_k\nabla f(y^{k+1})_i-(x_i^*-\lambda_k h^*_i)\right|^2-\rho^2\lambda_k^2. \label{eq:nice}
  \label{eq:mid}
  \end{eqnarray}
  Inequality (\ref{eq:nice}) follows because of the following fact: 
  \begin{eqnarray}
  a\geq b\geq 0\implies(a-b)^2\leq a^2-b^2\label{hey}
  \end{eqnarray}
  which applies because 
  \begin{eqnarray}
  \label{gh}
  |(y_i^{k+1}-\lambda_k\nabla f(y^{k+1})_i)-(x_i^*-\lambda_k h^*_i)|&\geq |(x_i^*-\lambda_k h^*_i)|
  \label{gg}
  \geq\rho\lambda_k>0
  \end{eqnarray}
  where we have used (\ref{eq:sgn}) and Lemma \ref{lemmaeq12} to prove (\ref{gg}).
  
  Now define for $k\in\mathbb{N}$,
  $$
  \label{Pdef}
  \calP_k\triangleq \{i\in E : \sgn(y_i^{k}-\lambda_{k-1}\nabla f(y^{k})_i)\neq -h_i^*/\rho\}
  $$ and recall the standard notation $|\calP_k|$ for the number of elements in $\calP_k$. 
  For all $k\in\mathbb{N}$:
  \begin{eqnarray}
  \nonumber\|x^{k+1}-x^*\|^2
  &=&
  \sum_{j\in\calP_{k+1}}|x_j^{k+1}-x^*_j|^2
  +\sum_{j\in [n]\backslash\calP_{k+1}}|x^{k+1}_j-x^*_j|^2 \\
  &\leq &
   \sum_{j\in\calP_{k+1}}\left\{
   \left|y_j^{k+1}-\lambda_k\nabla f(y^{k+1})_j-(x_j^*-\lambda_k h^*_j)\right|^2-\rho^2\lambda_k^2\right\}
  \nonumber\\&&\label{eq:l1lsnonx}
  +\sum_{j\in[n]\backslash\calP_{k+1}}\left|y_j^{k+1}-\lambda_k\nabla f(y^{k+1})_j-(x_j^*-\lambda_k h^*_j)\right|^2
  \\\nonumber
  &=&\|y^{k+1}-\lambda_k\nabla f(y^{k+1}) - (x^*-\lambda_k h^*)\|^2-\rho^2\lambda_k^2|\calP_{k+1}|\\                 
  &\leq&\|y^{k+1}-x^*\|^2-\rho^2\lambda_1^2|\calP_{k+1}|\label{eq:nonexpand}\label{eq:thebiggyv0}
  \label{eq:thebiggyv1}.
  \end{eqnarray}
  Inequality (\ref{eq:l1lsnonx}) follows from (\ref{eq:mid}) and the elementwise nonexpansiveness of $S_{\rho\lambda_k}$ (i.e. Lemma \ref{lemma:Sv}(\ref{eq:lemmaSv1})). To deduce (\ref{eq:nonexpand}), we used the fact that $I-\lambda\nabla f$ is nonexpansive for $0<\lambda<2/L$ \cite[Pro. 4.33]{bauschke2011convex}, and $\{\lambda_k\}_{k\in\mathbb{N}}$ is nondecreasing. Now (\ref{eq:thebiggyv1}) is in the form of (\ref{eq:lem_recurs}) of Lemma \ref{lem_recurse} with $x=x^*$ and $N_k = \rho^2\lambda_1^2 |\mathcal{P}_{k}|$. Since we assumed $\rho>0$ and $\lambda_1>0$ it follows that $\sum_{k\in\mathbb{N}} |\calP_k| <\infty$ for either parameter choice 1 or 2.  This implies $|\calP_k|$ is nonzero for only finitely many iterations, thus (\ref{thm:finite:eq2}) is proved. 
  
  \subsubsection*{Proof of (\ref{thm:finite:eq1}) of Theorem \ref{thm:finite}} 
  For $E^c$ nonempty, define
  $
  \label{omegaDef}
  \omega\triangleq\min\{\rho-|h^*_i|:i\in E^c\}\in(0,\rho].
  $
  If $E^c$ is empty, (\ref{thm:finite:eq1}) is trivially true, therefore assume $E^c$ is nonempty and note that 
  \begin{eqnarray} 
  \omega\lambda_k =\min\{\rho\lambda_k-\lambda_k |h^*_i|:i\in E^c\}>0.
  \label{mm}
  \end{eqnarray}
  Consider $i\in E^c$ (which implies $i\notin\supp(x^*)$). If $x_i^{k+1}\neq0$, then Lemma \ref{lemma:Sv} (\ref{eq:lemmaSv3}), (\ref{xupdate_special_ifbs}) and optimality condition (\ref{optcond_fixed}) imply
  \begin{eqnarray}
  \nonumber
  |x_i^{k+1}|^2
  &=& 
  |S_{\rho\lambda_{k}}\lbr y_i^{k+1}-\lambda_k\nabla f(y^{k+1})_i\rbr 
  - S_{\rho\lambda_{k}}\lbr - \lambda_k h_i^*\rbr|^2
  \\\nonumber
  &\leq&
  \left[
  |y_i^{k+1}-\lambda_k\nabla f(y^{k+1})_i+\lambda_k h^*_i|
  -\lbr\rho\lambda_k - \lambda_k | h_i^*|\rbr
  \right]^2
  \\\label{mnn}
  &\leq&
  |y_i^{k+1}-\lambda_k\nabla f(y^{k+1})_i+\lambda_k h^*_i|^2
  -\lbr\rho\lambda_k - \lambda_k |h_i^*|\rbr^2
  \\\label{mmn}
  &\leq&
  |y_i^{k+1}-\lambda_k\nabla f(y^{k+1})_i+\lambda_k h^*_i|^2
  - \omega^2\lambda^2_k.
  \end{eqnarray}
  To derive (\ref{mmn}) we used (\ref{mm}). To derive (\ref{mnn}) we used (\ref{hey}) which applies because 
  \begin{eqnarray}
  \label{mimo}
  |y_i^{k+1}-\lambda_k\nabla f(y^{k+1})_i+\lambda_k h^*_i|&\geq &|y_i^{k+1}-\lambda_k\nabla f(y^{k+1})_i|-\lambda_k|h_i^*|
  \\\label{meme}
  &>&\rho\lambda_k - \lambda_k| h_i^*|
  \end{eqnarray}
  which is greater than $0$ by (\ref{mm}). Note that (\ref{mimo}) follows from the identity: 
  $$|a+b|\geq|a|-|b|,\quad\forall\, a,b\in\mathbb{R}$$
  and (\ref{meme}) follows from the fact that $0\neq x_i^{k+1}=S_{\rho\lambda_k}(y_i^{k+1}-\lambda_k\nabla f(y^{k+1}))_i$. 
  Analogously to the definition of $\calP_k$, define for all $k\in\mathbb{N}$,
  $$
  \mathcal{Q}_k \triangleq\{i\in E^c: x_i^k\neq 0\}.
  $$
  Thus for all $k\in\mathbb{N}$,
  \begin{eqnarray}
  \|x^{k+1}-x^*\|^2
  &=&
  \sum_{j\in[n]\backslash\calQ_{k+1}}|x_j^{k+1} - x^*_j|^2
  +
  \sum_{j\in\calQ_{k+1}}|x_j^{k+1}|^2
  \nonumber\\
  &\leq&
  \sum_{j\in[n]\backslash\calQ_{k+1}}|y_j^{k+1}-\lambda_k\nabla f(y^{k+1})_j - (x^*_j-\lambda_k h_j^*)|^2
  \nonumber\\
  &&
  +
  \sum_{j\in\calQ_{k+1}}\left\{|y_j^{k+1}-\lambda_k\nabla f(y^{k+1})_j +\lambda_k h_j^*|^2-\omega^2\lambda_k^2\right\}
  \nonumber\\
  &\leq &
  \|y^{k+1}-x^*\|^2
  -\omega^2\lambda_1^2|\calQ_{k+1}|.
  \label{eq:thebiggy2}\nonumber
  \end{eqnarray}
  This recursion is in the form of (\ref{eq:lem_recurs}) in Lemma \ref{lem_recurse} with $x=x^*$ and $N_k = \omega^2\lambda_1^2|\calQ_{k}|$. Since $\omega$ and $\lambda_1$ are both greater than $0$ we have $\sum_{k\in\mathbb{N}} |\calQ_{k}|<\infty$. Thus $\calQ_k$ is nonempty for only finitely many iterations. Note that by (\ref{IFBS_yupdate}), if $x^k_i$ and $x^{k-1}_i$ are equal to $0$, then $y_i^{k+1}=0$. Thus (\ref{thm:finite:eq1}) is proved. \qed
  
  
  
  
  \section{Proof of Theorem \ref{cor:q_lin}}\label{Proof_q_lin}
  We first prove (\ref{eq:proj_HB2}). From Theorem \ref{thm:finite}, there exists $K>0$ such that for all $k>K$, (\ref{thm:finite:eq2}) and (\ref{thm:finite:eq1}) hold for either parameter choice 1 or 2. For $i\in E$, $k>K$, we calculate the quantity
  \begin{eqnarray}
  z_i^{k+1}&\triangleq& y_i^{k+1}-\lambda_k\nabla\phi(y_E^{k+1})_i
  \nonumber\\\label{bla}
  &=& y_i^{k+1}-\lambda_k(-h^*_i+\nabla f(y^{k+1})_i)\\\nonumber
  &=& y_i^{k+1}-\lambda_k \nabla f(y^{k+1})_i+\rho\lambda_k\left(\frac{h_i^*}{\rho}\right)\\\label{omg}
  &=&\sgn\left(y_i^{k+1}-\lambda_k\nabla f(y^{k+1})_i\right)(|y_i^{k+1}-\lambda_k \nabla f(y^{k+1})_i|-\rho\lambda_k).
  \end{eqnarray}
  Equation (\ref{bla}) follows from $\supp(y^{k+1})\subseteq E$. Equation (\ref{omg}) follows from (\ref{thm:finite:eq2}). Therefore, for $i\in E$, $k>K$,
  \begin{eqnarray*}
  	x_i^{k+1}&=& S_{\rho\lambda_k}\left(y_i^{k+1}-\lambda_k \nabla f(y^{k+1})_i\right)=
  	\left\{
  	\begin{array}{lr}
  		z_i^{k+1}: -h_i^*z_i^{k+1}\geq 0 \\
  		0  : \text{else}
  	\end{array}
  	\right.
  \end{eqnarray*}
  which proves (\ref{eq:proj_HB2}). Now (\ref{thm:finite:eq2}) implies $\sgn(x_i^k)=-h_i^*/\rho$ for all $i\in E$, $k>K$, and $x_i^k\neq 0$. Therefore $-h_i^* x_i^k = \rho|x_i^k|$, for all $i\in E$, $k>K$.
  Therefore since $x_{E^c}^k=0$ for $k>K$, $-(h_E^*)^\top x_E^k = \rho\|x^k\|_1$, which implies $F(x^k)=\phi(x^k_E)$.\qed
  
  \section{Proofs of Sec.~\ref{sec:locallin2}}
  \label{proofOfTheorem55}
  \subsubsection*{Proof of Proposition \ref{lemmings}}
  Corollary \ref{cor_IFBS_weak} implies that $\lim_{k\to\infty} x^k\triangleq x^*$ exists and $x^*\in X^*$ for parameter choice 1. On the other hand Lemma \ref{lemCham} implies this is true for parameter choice 2. Theorem \ref{thm:finite} states that there exists a finite $K$ such that for $k>K$ (\ref{thm:finite:eq2}) holds for all $i\in E$, and recall that $\supp(x^*)\subseteq E$.
  Now since $x^k_i\to x^*_i\neq 0$ for all $i\in\supp(x^*)$, there exists some $K'>0$ such that for all $k>K'$, $x_i^k\neq 0$. Combining this with (\ref{xupdate_special_ifbs}), (\ref{soft_thresh_def}) and (\ref{thm:finite:eq2}) implies that for all $k>\max(K,K')$, and $i\in\supp(x^*)$,
  \begin{eqnarray*}
  	x_i^{k} &=& \sgn(y_i^{k} - \lambda_{k-1}\nabla f(y^k)_i)(|y_i^{k} - \lambda_{k-1}\nabla f(y^k)_i| - \lambda_{k-1}\rho)
  	\\
  	&=&
  	y_i^{k} - \lambda_{k-1}(\nabla f(y^k)_i - h_i^*).
  \end{eqnarray*}\qed
  
  \subsubsection*{Proof of Theorem \ref{thm:furtherconv}}
  Recall that Prob (\ref{prob:l1LS}) satisfies Assumption 1, therefore all conclusions of Lemma \ref{thm:constGrad} hold. Further recall that $X^*$ is nonempty for Prob.~(\ref{prob:l1LS}). Therefore $\lim_{k\to\infty} x^k\triangleq x^*$ exists and $x^*\in X^*$ by Corollary \ref{cor_IFBS_weak}. 
  Recall that $E=\{i:|h_i^*|=\rho\}$ and also by the strict complementarity assumption: $E=\supp(x^*)$.
  Proposition \ref{lemmings} proves that there exists $K>0$ such that for all $k>K$
  \begin{eqnarray}
  x^k_E = y^k_E -\frac{1}{L}(\nabla f(y^k)_E-h_E^*)
  =
  y^k_E - \frac{1}{L}((A^\top A y^k)_E - (A^\top A x^*)_E).\label{hh}
  \end{eqnarray}
  On the other hand Theorem \ref{thm:finite} proved that there exists $K'>0$ such that for all $k>K'$, $x_{E^c}^k=y_{E^c}^k=0$. Therefore for all $k>\max(K,K')\triangleq K''$ both conditions hold.
  Let $Q=(A_E^\top A_E)$ and $P$ be the orthogonal projector for the range space of $Q$. 
  
  We first consider the part of the error in the nullspace of $P$. Equation (\ref{hh}) implies
  \begin{eqnarray*}
  	(I-P)(x_E^k - x^*_E)=(I-P)(y_E^k - x_E^*),\quad\forall k>K''.
  \end{eqnarray*}
  Combining this with (\ref{IFBS_yupdate}) implies:
  $
  t^{k+1}=
  (1+\alpha)t^k
  -
  \alpha t^{k-1}
  $
  where $t^k = (I-P)(y_E^k-x_E^*)$ for all $k>K''$. This is a linear homogeneous recursion with solution:
  \begin{eqnarray*}
  	\tilde{t}^M_i = \tilde{t}^0_i + \frac{(\tilde{t}^1_i-\tilde{t}^0_i)(1-\alpha^M)}{1-\alpha},\quad\forall M\in\mathbb{N},\label{hardss}
  \end{eqnarray*}
  where $\tilde{t}^k = t^{k+\ceil{K''}}$. Now 
  $
  \lim_{M\to\infty}\tilde{t}_i^M = (\tilde{t}_i^1 - \alpha \tilde{t}_i^0)/(1-\alpha).
  $
  On the other hand,
  Thm. \ref{thm:fista_lyap} (\ref{state:weak}) implies $\tilde{t}^M_i\to 0$ as $M\to\infty$. Therefore either $\tilde{t}^M=0$ for all $M\in\mathbb{N}$ or $\tilde{t}_i^M=\alpha \tilde{t}_i^{M-1}$  for all $M$. Therefore $(I-P)(y^k_E-x_E^*)=(I-P)(x^k_E-x_E^*)$ converges to $0$ R-linearly with rate $\alpha$.
  
  Next we consider $P(x_E^k-x_E^*)$. Note that $Q$ is symmetric thus $\mathcal{R}(Q)=\mathcal{N}(Q)^{\perp}$. Thus, for all $k>K''$ 
  \begin{eqnarray}
 P( x_E^k-x_E^*) &=& P(y_E^k-x_E^*) - \frac{1}{L}PQ(y_E^k - x_E^*)
 \nonumber\\
  &=&
  P(y_E^k-x_E^*) - \frac{1}{L}Q P(y_E^k - x_E^*)\label{cc}.
  \end{eqnarray}
  Let $\hat{l}_E$ is the smallest eigenvalue of $Q$ restricted $\mathcal{R}(P)$. If $\hat{l}_E=0$, then $P$ is the all-zero matrix and $x^k$ converges to $x^*$ R-linearly with rate $\alpha$. Assume $\hat{l}_E>0$. Restating  (\ref{IFBS_yupdate}) we have for all $k\geq K''$: (\ref{cc}) holds and
  $$
  P y_E^{k+1} = P x_E^{k}+\alpha(P x_E^k-P x_E^{k-1}).\label{bright_recursion2}
  $$
  This is exactly the same recursion as studied in \cite[\S 4.2--4.3]{o2012adaptive} with respect to the sequences $\{P(x_E^k-x^*_E)\}$ and $\{P(y_E^k-x_E^*)\}$. Note that $\phi$ restricted to $\mathcal{R}(P)$ is a strongly-convex quadratic function.  By looking at the eigenvalues and eigenvectors of $Q$ restricted to $\mathcal{R}(P)$, one can see that Q-linear convergence of $P x_E^k$ is obtained and the rate $(1-\sqrt{\mu/L})^{1/2}$ is achieved by the choice: $\alpha = (1-\sqrt{\mu/L})/(1+\sqrt{\mu/L})$. We refer to \cite{o2012adaptive} for all the details.  Note that the rate of $x^k$ is the same as $P x_E^k$ since $x^k$ is zero outside $E$ for $k>K''$ and $(I-P)(x_E^k-x^*_E)$ has R-linear convergence to $0$ with rate $\alpha$, which is faster than the rate $(1-\sqrt{\mu/L})^{1/2}$. Finally the fact that $\phi$ is quadratic for this problem gives the objective function rates. 

\bibliographystyle{spmpsci} 
\bibliography{C:/Users/Pat/Dropbox/ResearchF14on/ICASSP14/refs}

\begin{thebibliography}{10}
\providecommand{\url}[1]{{#1}}
\providecommand{\urlprefix}{URL }
\expandafter\ifx\csname urlstyle\endcsname\relax
  \providecommand{\doi}[1]{DOI~\discretionary{}{}{}#1}\else
  \providecommand{\doi}{DOI~\discretionary{}{}{}\begingroup
  \urlstyle{rm}\Url}\fi

\bibitem{afonso2010fast}
Afonso, M.V., Bioucas-Dias, J.M., Figueiredo, M.A.: Fast image recovery using
  variable splitting and constrained optimization.
\newblock Image Processing, IEEE Transactions on \textbf{19}(9), 2345--2356
  (2010)

\bibitem{agarwal2012}
Agarwal, A., Negahban, S., Wainwright, M.J.: Fast global convergence of
  gradient methods for high-dimensional statistical recovery.
\newblock Ann. Statist. \textbf{40}(5), 2452--2482 (2012).
\newblock \doi{10.1214/12-AOS1032}.
\newblock \urlprefix\url{http://dx.doi.org/10.1214/12-AOS1032}

\bibitem{alvarez2000minimizing}
Alvarez, F.: {On the minimizing property of a second order dissipative system
  in Hilbert spaces}.
\newblock SIAM Journal on Control and Optimization \textbf{38}(4), 1102--1119
  (2000)

\bibitem{attouch2015_2}
Attouch, H., Chbani, Z., Peypouquet, J., Redont, P.: {Fast Convergence of
  Inertial Dynamics and Algorithms with Asymptotic Vanishing Damping}.
\newblock Preprint, available at
  \url{http://www.optimization-online.org/DB_FILE/2015/10/5179.pdf}  (2015)

\bibitem{attouch2015rate}
Attouch, H., Peypouquet, J.: {The rate of convergence of Nesterov's accelerated
  forward-backward method is actually $o (k^{-2})$}.
\newblock arXiv preprint arXiv:1510.08740  (2015)

\bibitem{attouch2014dynamical}
Attouch, H., Peypouquet, J., Redont, P.: {A dynamical approach to an inertial
  forward-backward algorithm for convex minimization}.
\newblock SIAM Journal on Optimization \textbf{24}(1), 232--256 (2014)

\bibitem{bach2011convex}
Bach, F., Jenatton, R., Mairal, J., Obozinski, G., et~al.: Convex optimization
  with sparsity-inducing norms.
\newblock Optimization for Machine Learning pp. 19--53 (2011)

\bibitem{bauschke2011convex}
Bauschke, H.H., Combettes, P.L.: Convex analysis and monotone operator theory
  in Hilbert spaces.
\newblock Springer Science \& Business Media (2011)

\bibitem{Beck:2009}
Beck, A., Teboulle, M.: {A fast iterative shrinkage-thresholding algorithm for
  linear inverse problems}.
\newblock SIAM J. Img. Sci. \textbf{2}(1), 183--202 (2009).
\newblock \doi{10.1137/080716542}.
\newblock \urlprefix\url{http://dx.doi.org/10.1137/080716542}

\bibitem{bertsekas1999nonlinear}
Bertsekas, D.P.: Nonlinear programming, 2nd edn.
\newblock {Athena Scientific} (1999)

\bibitem{boct2015inertial}
Bo{\c{t}}, R.I., Csetnek, E.R., Hendrich, C.: {Inertial Douglas--Rachford
  splitting for monotone inclusion problems}.
\newblock Applied Mathematics and Computation \textbf{256}, 472--487 (2015)

\bibitem{bredies2008linear}
Bredies, K., Lorenz, D.A.: {Linear convergence of iterative soft-thresholding}.
\newblock Journal of Fourier Analysis and Applications \textbf{14}(5-6),
  813--837 (2008)

\bibitem{burachik1999varepsilon}
Burachik, R.S., Svaiter, B.: $\varepsilon$-enlargements of maximal monotone
  operators in banach spaces.
\newblock Set-Valued Analysis \textbf{7}(2), 117--132 (1999)

\bibitem{cevher2014convex}
Cevher, V., Becker, S., Schmidt, M.: Convex optimization for big data:
  Scalable, randomized, and parallel algorithms for big data analytics.
\newblock Signal Processing Magazine, IEEE \textbf{31}(5), 32--43 (2014)

\bibitem{chambolle2015convergence}
Chambolle, A., Dossal, C.: On the convergence of the iterates of the “fast
  iterative shrinkage/thresholding algorithm”.
\newblock Journal of Optimization Theory and Applications pp. 1--15 (2015)

\bibitem{chambolle2014ergodic}
Chambolle, A., Pock, T.: {On the ergodic convergence rates of a first-order
  primal-dual algorithm.} (2015).
\newblock \urlprefix\url{https://hal.archives-ouvertes.fr/hal-01151629}.
\newblock To appear in Math. Programm. A

\bibitem{choi2010compressed}
Choi, K., Wang, J., Zhu, L., Suh, T.S., Boyd, S., Xing, L.: {Compressed sensing
  based cone-beam computed tomography reconstruction with a first-order
  method}.
\newblock Medical physics \textbf{37}(9), 5113--5125 (2010)

\bibitem{prox_signalProcessing}
Combettes, P.L., Pesquet, J.C.: {Proximal splitting methods in signal
  processing}.
\newblock In: Fixed-Point Algorithms for Inverse Problems in Science and
  Engineering, pp. 185--212. Springer (2011)

\bibitem{combettes2005signal}
Combettes, P.L., Wajs, V.R.: Signal recovery by proximal forward-backward
  splitting.
\newblock Multiscale Modeling \& Simulation \textbf{4}(4), 1168--1200 (2005)

\bibitem{condat2013primal}
Condat, L.: {A primal--dual splitting method for convex optimization involving
  Lipschitzian, proximable and linear composite terms}.
\newblock Journal of Optimization Theory and Applications \textbf{158}(2),
  460--479 (2013)

\bibitem{davis2014convergence}
Davis, D., Yin, W.: Convergence rate analysis of several splitting schemes.
\newblock Tech. rep., UCLA CAM Report 14-51 (2014)

\bibitem{eckstein2012augmented}
Eckstein, J., Yao, W.: {Augmented Lagrangian and alternating direction methods
  for convex optimization: A tutorial and some illustrative computational
  results}.
\newblock RUTCOR Research Reports \textbf{32} (2012)

\bibitem{friedman2010regularization}
Friedman, J., Hastie, T., Tibshirani, R.: Regularization paths for generalized
  linear models via coordinate descent.
\newblock Journal of statistical software \textbf{33}(1), 1 (2010)

\bibitem{guler1992new}
G{\"u}ler, O.: New proximal point algorithms for convex minimization.
\newblock SIAM Journal on Optimization \textbf{2}(4), 649--664 (1992)

\bibitem{Hale:2008}
Hale, E.T., Yin, W., Zhang, Y.: {Fixed-point continuation for
  $\ell_1$-minimization: methodology and convergence}.
\newblock SIAM J. on Optimization \textbf{19}(3), 1107--1130 (2008).
\newblock \doi{10.1137/070698920}

\bibitem{hare2004identifying}
Hare, W., Lewis, A.S.: Identifying active constraints via partial smoothness
  and prox-regularity.
\newblock Journal of Convex Analysis \textbf{11}(2), 251--266 (2004)

\bibitem{kim2007interior}
Kim, S.J., Koh, K., Lustig, M., Boyd, S., Gorinevsky, D.: {An interior-point
  method for large-scale $\ell_1$-regularized least squares}.
\newblock Selected Topics in Signal Processing, IEEE Journal of \textbf{1}(4),
  606--617 (2007)

\bibitem{liang2014local}
Liang, J., Fadili, J., Peyr{\'e}, G.: {Local linear convergence of
  forward--backward under partial smoothness}.
\newblock In: Advances in Neural Information Processing Systems, pp. 1970--1978
  (2014)

\bibitem{liang2015activity}
Liang, J., Fadili, J., Peyr{\'e}, G.: Activity identification and local linear
  convergence of inertial forward-backward splitting.
\newblock arXiv preprint arXiv:1503.03703  (2015)

\bibitem{liang2014convergence}
Liang, J., Fadili, J., Peyr{\'e}, G.: Convergence rates with inexact
  nonexpansive operators.
\newblock Accepted Mathematical Programming  (2015)

\bibitem{lorenz2014inertial}
Lorenz, D.A., Pock, T.: An inertial forward-backward algorithm for monotone
  inclusions.
\newblock Journal of Mathematical Imaging and Vision pp. 1--15 (2014).
\newblock \doi{10.1007/s10851-014-0523-2}.
\newblock \urlprefix\url{http://dx.doi.org/10.1007/s10851-014-0523-2}

\bibitem{lustig2007sparse}
Lustig, M., Donoho, D., Pauly, J.M.: {Sparse MRI: The application of compressed
  sensing for rapid MR imaging}.
\newblock Magnetic resonance in medicine \textbf{58}(6), 1182--1195 (2007)

\bibitem{mainge2008convergence}
Maing{\'e}, P.E.: {Convergence theorems for inertial KM-type algorithms}.
\newblock Journal of Computational and Applied Mathematics \textbf{219}(1),
  223--236 (2008)

\bibitem{monteiro2012adaptive}
Monteiro, R.D., Ortiz, C., Svaiter, B.F.: An adaptive accelerated first-order
  method for convex optimization.
\newblock Tech. rep., School of Industrial and Systems Engineering, Georgia
  Institute of Technology, Atlanta, GA (2012)

\bibitem{MoudafiOliny}
Moudafi, A., Oliny, M.: Convergence of a splitting inertial proximal method for
  monotone operators.
\newblock Journal of Computational and Applied Mathematics \textbf{155}(2), 447
  -- 454 (2003).
\newblock \doi{http://dx.doi.org/10.1016/S0377-0427(02)00906-8}

\bibitem{nesterov1983method}
Nesterov, Y.: {A method of solving a convex programming problem with
  convergence rate $O(1/k^2)$}.
\newblock In: Soviet Mathematics Doklady, vol.~27, pp. 372--376 (1983)

\bibitem{nesterov2004introductory}
Nesterov, Y.: {Introductory lectures on convex optimization: a basic course}.
\newblock Springer (2004)

\bibitem{nesterov2013gradient}
Nesterov, Y.: Gradient methods for minimizing composite functions.
\newblock Mathematical Programming \textbf{140}(1), 125--161 (2013)

\bibitem{o2012adaptive}
O'Donoghue, B., Cand{\`e}s, E.: {Adaptive restart for accelerated gradient
  schemes}.
\newblock Foundations of Computational Mathematics pp. 1--18 (2012)

\bibitem{opial1967weak}
Opial, Z.: {Weak convergence of the sequence of successive approximations for
  nonexpansive mappings}.
\newblock Bulletin of the American Mathematical Society \textbf{73}(4),
  591--597 (1967)

\bibitem{passty1979ergodic}
Passty, G.B.: {Ergodic convergence to a zero of the sum of monotone operators
  in Hilbert space}.
\newblock Journal of Mathematical Analysis and Applications \textbf{72}(2),
  383--390 (1979)

\bibitem{polyak1964some}
Polyak, B.T.: {Some methods of speeding up the convergence of iteration
  methods}.
\newblock USSR Computational Mathematics and Mathematical Physics
  \textbf{4}(5), 1--17 (1964)

\bibitem{PolyakIntro}
Polyak, B.T.: {Introduction to Optimization}.
\newblock Optimization Software Inc. (1987)

\bibitem{raguet2013generalized}
Raguet, H., Fadili, J., Peyr{\'e}, G.: A generalized forward-backward
  splitting.
\newblock SIAM Journal on Imaging Sciences \textbf{6}(3), 1199--1226 (2013)

\bibitem{shalev2011stochastic}
Shalev-Shwartz, S., Tewari, A.: Stochastic methods for $\ell_1$-regularized
  loss minimization.
\newblock The Journal of Machine Learning Research \textbf{12}, 1865--1892
  (2011)

\bibitem{su2014differential}
Su, W., Boyd, S., Cand{\`e}s, E.: {A differential equation for modeling
  Nesterov's accelerated gradient method: theory and insights}.
\newblock In: Advances in Neural Information Processing Systems, pp. 2510--2518
  (2014)

\bibitem{tao2015local}
Tao, S., Boley, D., Zhang, S.: {Local Linear Convergence of ISTA and FISTA on
  the LASSO Problem}.
\newblock Tech. rep., University of Minnesota (2015)

\bibitem{tibshirani2013lasso}
Tibshirani, R.J.: The lasso problem and uniqueness.
\newblock Electronic Journal of Statistics \textbf{7}, 1456--1490 (2013)

\bibitem{tseng2008accelerated}
Tseng, P.: {On accelerated proximal gradient methods for convex-concave
  optimization. Submitted to SIAM J. Opt.}  (2008)

\bibitem{wen2012convergence}
Wen, Z., Yin, W., Zhang, H., Goldfarb, D.: On the convergence of an active-set
  method for $\ell_1$ minimization.
\newblock Optimization Methods and Software \textbf{27}(6), 1127--1146 (2012)

\bibitem{zhang2015necessary}
Zhang, H., Yin, W., Cheng, L.: Necessary and sufficient conditions of solution
  uniqueness in 1-norm minimization.
\newblock Journal of Optimization Theory and Applications \textbf{164}(1),
  109--122 (2015)

\end{thebibliography}

\end{document}